\documentclass[11pt,reqno]{amsart}
\usepackage[width=150mm,height=220mm,
centering,
]{geometry}
\usepackage{mathrsfs}
\usepackage{amsxtra}
\usepackage{amsmath,amsfonts,amssymb}
\usepackage{upgreek}
\usepackage{xcolor}
\usepackage[numbers,sort&compress]{natbib}
\usepackage{csquotes}
\numberwithin{equation}{section}
\allowdisplaybreaks
\ifcsname proof\endcsname
  \else
\newcommand{\proof}[1][Proof]{\textbf{#1}.\,\,}
  \def\endproof{\hfill $\Box$}
\fi

\newtheorem{Theorem}{Theorem}[section]

\newtheorem{Lemma}[Theorem]{Lemma}

\newtheorem{Remark}[Theorem]{Remark}

\def\rmd{\mathrm{d}}
\begin{document}

\parindent 2em
\title{On the rate of convergence for Landau type Schr\"odinger Operators}
\author{Yucheng Pan}
\address{School of Mathematical Sciences and LPMC,
  Nankai University, Tianjin,  China}
\email{panyucheng@mail.nankai.edu.cn}

\author{Wenchang Sun}
\address{School of Mathematical Sciences and LPMC,
  Nankai University, Tianjin,  China}
\email{sunwch@nankai.edu.cn}

\thanks{Corresponding author: Wenchang Sun.}
\thanks{This work was supported by the
National Natural Science Foundation of China (No. 12171250, U21A20426 and 12271267).}
\date{}

\begin{abstract}
We study the pointwise convergence of Landau type Schr\"odinger operators within the fractional Sobolev space $W^{s,p}(\mathbb R)$.
Our results extend those established by Bailey (Rev. Mat. Iberoam., 29 (2): 531–546, 2013)
and Yuan, Zhao and Zheng (Nonlinear Anal., 208: Paper No. 112312, 28, 2021).
Furthermore, we also analyze the convergence rate of Landau type Schr\"odinger operators along curves
and derive a sharp result for the case of convergence along vertical lines.
\end{abstract}

\subjclass[2020]{46E35}
\keywords{Landau type Schr\"odinger operator, restricted curve, convergence rate}

\maketitle

\section{Introduction and Main Results}

The Schr\"odinger equation stands as one of the cornerstones of quantum mechanics,
revealing the fundamental principles that govern the behavior of matter in the microscopic realm.
It serves as an indispensable tool for addressing almost all non-relativistic issues within atomic physics.
Since
this equation was introduced
by Schr\"odinger \cite{schrodinger1926}
 in 1926,
it has found extensive applications in various fields,
including solid-state physics, nuclear physics and chemistry.

It is well known that when the initial datum is sufficiently well behaved
(for instance, being a Schwartz function), the solution to the Schr\"odinger equation can be written as
\begin{equation}\label{eq:s1t0}
  e^{it(-\Delta)^{a/2}}f(x)=\frac{1}{(2\pi)^n}\int_{\mathbb R^n}e^{ix\cdot \xi}e^{i g(t)|\xi|^a}\hat{f}(\xi)\rmd \xi ,
\end{equation}
where $a>0$ and $g(t)=t$. The operator $e^{it(-\Delta)^{a/2}}$ defined by (\ref{eq:s1t0}) is referred to as the fractional Schr\"odinger operator.

The research on the pointwise convergence problem focuses on determining the optimal regularity index $s$ for which
\begin{equation}\label{eq:s1t9}
  \lim_{t\rightarrow 0}e^{it(-\Delta)^{a/2}}f(x)=f(x), \quad a.e. \ x\in \mathbb R^n ,
\end{equation}
whenever $f$ belongs to the Sobolev space $H^s (\mathbb R^n)$.
For the case $a=2$ and $n=1$, Carleson \cite{Carleson1980} first considered this problem and proved that the convergence holds when $s\geq 1/4$.
Later on, Dahlberg and Kenig \cite{Dahlberg1982} constructed a counterexample showing that the range of regularity index given by Carleson is sharp.
In recent years, driven by advancements in restriction estimates within Fourier analysis,
the technique of polynomial partitioning and decoupling has gradually gained recognition.
By decoupling, Du,  Guth, and  Li \cite{Guth2017} showed that $s>1/2-1/(2n+2)$ is sufficient for (\ref{eq:s1t9}) to hold for the dimension $n=2$
and Du and Zhang \cite{DuZhang2019} proved that this range also holds for general $n\geq 3$.
In fact, their results are sharp up to the endpoint since Bourgain \cite{Bourgain2016} proved that the convergence does not hold for $s<1/2-1/(2n+2)$.
See also works by Sj\"olin \cite{Sjolin1987} and Walther \cite{Walther2001} on this topic.

As an extension of the convergence problem associated with fractional Schr\"odinger operators,
many authors have also focused on the convergence properties of solutions to other PDEs arising in physical contexts.
For example, if we take $g(t)=t+it^{\gamma}$ and $n=1$ in (\ref{eq:s1t0}), then we get the Landau type Schr\"odinger operator defined by
\begin{equation}\label{eq:s1t10}
  P_{a,\gamma}^t f(x):=\frac{1}{2\pi}\int_{\mathbb R}\hat{f}(\xi)e^{ix \xi}e^{it|\xi|^a}e^{-t^\gamma |\xi|^a}\rmd \xi .
\end{equation}
Note that when $\gamma=1$ and $a=2$, $P_{a,\gamma}^{t} f$ corresponds to the solution of the Ginzburg-Landau equation \cite{Duan1994}.
This equation was initially introduced by Ginzburg and Landau in the 1950s during their study of phase transition phenomena in superconductors.
Presently, it is widely utilized to model a variety of phenomena,
including superconductivity, chemical turbulence, and various types of fluid flows (see \cite{Cazenave2014} for example).

Sj\"olin \cite{sjolin2009maximal}   initially investigated
the pointwise convergence of Landau type Schr\"odinger operators for the case $a=2$.
As an extension of the above result, for any $\gamma>1$ and $a>1$,
Bailey \cite{bailey2013boundedness} showed that the convergence holds for all $f\in H^s (\mathbb R)$ where $s>\min\{a(1-1/\gamma)/4,1/4\}$.
For the case $\gamma>1$ and $0<a\leq 1$,
Yuan,   Zhao  and   Zheng \cite{yuan2021dimension} (in the case $0<a\leq 1$) and Niu and Xue \cite{Niu2020} (in the case $0<a<1$) proved that
the convergence is still vaild when $s>a(1-1/\gamma)/4$ with $0<a<1$ or $s>(1-1/\gamma)/2$ with $a=1$.
Note that the ranges of the regularity index $s$ provided in \cite{bailey2013boundedness,yuan2021dimension,Niu2020} are all sharp up to the endpoint.
Moreover, for the case $0<\gamma \leq 1$ and $0<a\leq 1$,
Yuan,   Zhao  and   Zheng  \cite{yuan2021dimension} showed that the convergence holds for all $f\in L^p (\mathbb R)$, $1\leq p<\infty$.
Related works also include the convergence along curves by Niu and Xue \cite{Niu2020curve},
the non-tangential convergence by  Yuan,   Zhao  and   Zheng \cite{yuan2021pointwise}.
For the discrete case, the authors in \cite{Pan2025} studied the pointwise convergence of Landau type Schr\"odinger operators along sequences.

We first make the following observation.
It is obvious that (\ref{eq:s1t10}) makes sense for all $f\in L^p (\mathbb R)$, $1\leq p <2$.
Consequently, it is natural to ask whether the convergence still holds for $f\in L^p (\mathbb R)$ when $1\leq p <2$ and $\gamma >1$.
To this end, we present the following Theorem \ref{thm:s1t1} which provides a negative result regarding this question.

\begin{Theorem}\label{thm:s1t1}
Let $a>0$, $\gamma >1$ and $1\leq p<2$.
Then there exists a set $E\subset \mathbb R$ of positive measure and some $f_0 \in L^p (\mathbb R)$ such that
\begin{equation}\label{eq:s1t1}
 \lim_{t\rightarrow 0}P_{a,\gamma}^t f_0 (x)\neq f_0 (x),\quad \forall x\in E.
\end{equation}
\end{Theorem}

Based on Theorem \ref{thm:s1t1}, we proceed to study the pointwise convergence of
Landau type Schr\"odinger operators within the fractional Sobolev space $W^{s,p}\mathbb (\mathbb R)$ and obtain the following result.
It is easy to see that when $p=2$,
Theorem \ref{thm:s1t2} coincides with the results of Bailey \cite{bailey2013boundedness} and Yuan,   Zhao  and   Zheng \cite{yuan2021dimension}.

\begin{Theorem}\label{thm:s1t2}
Let $s>0$, $\gamma>1$ and $1<p\leq 2$.
If $0<a<1$ with $s>(2-p)/(2p) +a(1-1/\gamma)/4$,
and $a=1$ with $s>1/p-1/(2\gamma)$,
and $a>1$ with $s>(2-p)/(2p) +\min\{a(1-1/\gamma)/4,1/4\}$,
then for any $f\in W^{s,p}\mathbb (\mathbb R)$, we have
\begin{equation}\label{eq:s1t2}
 \lim_{t\rightarrow 0} P_{a,\gamma}^t f(x)=f(x), \quad a.e. \ x\in \mathbb R .
\end{equation}
\end{Theorem}

In recent years, many authors have also focused on the convergence rate for Schr\"odinger operators along curves.
Let $0<\beta\leq 1$ and $R>0$ be some constants.
Suppose that $x_0\in \mathbb R^n$ and the function $\Gamma : \mathbb R^n \times [0,1]\rightarrow \mathbb R^n$ satisfy
\begin{equation}\label{eq:s1t3}
 \Gamma (x,0)=x \ \mbox{ and } \
|\Gamma(x,t)-\Gamma(x,t')|\lesssim |t-t'|^\beta
\end{equation}
uniformly for $x\in  B(x_0,R)$ and $t,t'\in [0,1]$.
Additionally, assuming that the corresponding maximal estimates for Schr\"odinger operators have been established.
The convergence rate problem is to determine the optimal convergence rate $I(\beta,\delta,a)$
for all $f\in H^{s+\delta} (\mathbb R^n)$, $\delta \geq 0$, such that
\[
  e^{it(-\Delta)^{a/2}}(f)(\Gamma (x,t))-f(x)=o(t^{I(\beta,\delta,a)}),\quad a.e. \ x\in B(x_0 ,R)\quad as \quad t\rightarrow 0^+ .
\]
For the case of convergence along vertical lines in one spatial dimension,
 Cao,  Fan  and   Wang \cite{Cao2018} proved that when $a>1$, $s\geq 1/4$, and $0\leq \delta <a$, the convergence rate is $I=\delta /a$.
Different from the convergence rate of fractional Schr\"odinger operators discussed in \cite{Cao2018},
Li,  Niu and  Xue \cite{Li2024} investigated the rate of convergence along vertical lines for Landau type Schr\"odinger operators.
For higher dimensions,
Li and Wang \cite{Li2021JFA} studied the convergence rate along curves for a class of Schr\"odinger operators with polynomial growth.
Recently, by applying an effective time-frequency decomposition method, Li and Wang \cite{Li2021arxiv} have improved the above result to the optimal.

In this paper, we discuss the convergence rate for the Landau type Schr\"odinger operator along curves.
To be specific, we have the following result.
\begin{Theorem}\label{thm:s1t3}
Let $\gamma >0$, $a>0$ and $p\geq 1$.
Assume that the fuction $\Gamma :\mathbb R \times [0,1]\rightarrow \mathbb R$ satisfies (\ref{eq:s1t3}) for some $0<\beta\leq 1$.
If there exists $s_0 \geq 0$ such that for each $s>s_0$,
\begin{equation}\label{eq:s1t4}
\biggr\|\sup_{0<t<1}|P_{a,\gamma}^t (f)(\Gamma(x,t))|\biggr\|_{L^p (B(x_0 ,R))}\lesssim \|f\|_{H^s (\mathbb R)},
\end{equation}
then we have\\
$(1)$ when $a\geq \min\{1,\gamma\}$ and $\min\{1,\gamma\}/a \leq \beta \leq 1$, for each $s>s_0$ and $f\in H^{s+\delta}(\mathbb R)$, it holds,
\begin{equation}\label{eq:s1t5}
P_{a,\gamma}^t (f)(\Gamma(x,t))-f(x)=o(t^{h}),\quad a.e. \ x\in B(x_0 ,R)\quad as \quad t\rightarrow 0^+ ,
\end{equation}
whenever $(\delta,h)\in D_1 :=\{(\delta,h): \delta\geq 0, h\geq 0, h\leq \delta\cdot \min\{1,\gamma\}/a, h< \min\{\beta,\gamma\}\}$;\\
$(2)$ when $a\geq \min\{1,\gamma\}$ with $0<\beta<\min\{1,\gamma\}/a$, or $a< \min\{1,\gamma\}$ with $0<\beta\leq 1$,
for each $s>s_0$ and $f\in H^{s+\delta}(\mathbb R)$, it holds,
\begin{equation}\label{eq:s1t6}
P_{a,\gamma}^t (f)(\Gamma(x,t))-f(x)=o(t^{h}),\quad a.e. \ x\in B(x_0 ,R)\quad as \quad t\rightarrow 0^+ ,
\end{equation}
whenever $(\delta,h)\in D_2 :=\{(\delta,h): \delta\geq 0, h\geq 0, h\leq \beta \delta, h<\min\{\beta,\gamma\}\}$.
\end{Theorem}

According to Theorem \ref{thm:s1t3}, when $a\geq \min\{1,\gamma\}$,
the convergence rate for Landau type Schr\"odinger operators along vertical lines is faster than $t^{\delta\cdot \min\{1,\gamma\}/a}$,
while for $a< \min\{1,\gamma\}$, the convergence rate is faster than $t^{\delta}$.
In fact, we obtain a stronger result for the case $a< \min\{1,\gamma\}$, see Theorem \ref{thm:s1t4} below.
\begin{Theorem}\label{thm:s1t4}
Let $\gamma >0$, $a>0$ and $p\geq 1$. If there exists $s_0 \geq 0$ such that for each $s>s_0$,
\begin{equation}\label{eq:s1t17}
\biggr\|\sup_{0<t<1}|P_{a,\gamma}^t f|\biggr\|_{L^p (B(0,1))}\lesssim \|f\|_{H^s (\mathbb R)},
\end{equation}
then for all $f\in H^{s+\delta}(\mathbb R)$, $0\leq \delta <a$,
\begin{equation}\label{eq:s1t18}
P_{a,\gamma}^t f(x)-f(x)=o(t^{\delta\min\{1,\gamma\}/a}),\quad a.e. \ x\in \mathbb R \quad as \quad t\rightarrow 0^+ .
\end{equation}
\end{Theorem}
For the existence of $s_0$ in Theorem \ref{thm:s1t4}, we refer to
\cite[Theorem 1.2]{bailey2013boundedness} and \cite[Theorem 1.1]{yuan2021dimension}.
We further note that the proofs of Theorem \ref{thm:s1t3} and Theorem \ref{thm:s1t4} are both applicable to higher dimensions.
Since the proof of Theorem \ref{thm:s1t4} is very similar to that of Theorem \ref{thm:s1t3} $(1)$, we put it in the appendix.
Moreover, we also discuss the sharpness of Theorem \ref{thm:s1t4} in Section $2.4$.

\subsection*{Definitions and Notations}

Recall that for $0<s<1$ and $1\leq p<\infty$, the fractional Sobolev space $W^{s,p}(\mathbb R)$ endowed with the norm
\[
  \|u\|_{W^{s,p}(\mathbb R)}
  :=\biggr(\|u\|_{L^p (\mathbb R)}^p +\iint_{\mathbb R \times \mathbb R}\frac{|u(x)-u(y)|^p}{|x-y|^{sp+1}}\rmd x\rmd y \biggr)^{1/p}.
\]
For the case $s=k\geq 1$ be an integer, the norm $\|\cdot\|_{W^{k,p}(\mathbb R)}$ is defined by
\[
  \|u\|_{W^{k,p}(\mathbb R)}:=\biggr(\|u\|_{L^p (\mathbb R)}^p +\sum_{1\leq \alpha\leq k}\|D^{\alpha} u\|_{L^p (\mathbb R)}^p \biggr)^{1/p},
\]
where $D^{\alpha} u$ is the distributional derivative of $u$.
When $s>1$ not be an integer, we decompose $s$ as $s=m+\tau$, where $m\geq 1$ be an integer and $\tau \in (0,1)$. Then
\[
   \|u\|_{W^{s,p}(\mathbb R)}
  :=\biggr(\|u\|_{W^{m,p} (\mathbb R)}^p +\iint_{\mathbb R \times \mathbb R}\frac{|D^m u(x)-D^m u(y)|^p}{|x-y|^{\tau p+1}}\rmd x\rmd y \biggr)^{1/p}.
\]

Throughout this article, \enquote{$\lceil a\rceil$} represents the ceiling function of $a$.
By $M\lesssim N$ we mean that $M\leq cN$ for some constant $c$ independent of the parameters related to $M$ and $N$.
Besides, $M\approx N$ means that $M$ and $N$ are comparable, that is, $M\lesssim N$ and $N\lesssim M$.

\section{Proofs of Main Results}

\subsection{Proof of Theorem \ref{thm:s1t1}}

To prove Theorem \ref{thm:s1t1}, we first establish the connection between the almost everywhere convergence result
and the weak boundedness of the maximal operator $P_{a,\gamma}^{*}$ defined by
\[
  P_{a,\gamma}^{*}f(x)=\sup_{0<t<1}|P_{a,\gamma}^{t}f(x)|.
\]
To be specific, we have the following result.
\begin{Theorem}\label{thm:s2t1}
Let $a>0$, $\gamma >1$ and $1\leq p<2$. Suppose that the pointwise convergence (\ref{eq:s1t2}) holds for any $f\in L^p (\mathbb R)$.
Then for any $\varepsilon >0$, there exists a set $E_{\varepsilon}\subset [0,1]$ with $|E_{\varepsilon}|>1-\varepsilon$ such that
\begin{equation}\label{eq:s2t2}
|\{x\in E_{\varepsilon}:P_{a,\gamma}^{*}f(x)> \lambda\}|\leq C_{\varepsilon}\lambda^{-p}\|f\|_p^p
\end{equation}
holds for all $\lambda >0$ and $f\in L^p (\mathbb R)$.
\end{Theorem}
In the proof of Theorem \ref{thm:s2t1}, we need the following lemma given by Nikishin \cite{Nikishin1972}.
\begin{Lemma}[{\cite[Theorem 2]{Nikishin1972}}]\label{lem:s2t1}
Assume $X$ is a space with $\sigma$-finite measure, $D$ is an $N$-dimensional region with $|D|<\infty$ and $S(D)$ denotes the set of the measurable functions on the region $D$. Assume also that $L^p (X)$ is separable. Let $G$ be a bounded hyperlinear operator from $L^p (X)$ into $S(D)$, which means that
\[
  G(f+g)\leq Gf+ Gg,\quad \forall f,g\in L^p (X).
\]
Then for an arbitrary $\varepsilon >0$ there exists a set $E_{\varepsilon}\subset D$ with $|E_{\varepsilon}|\geq |D|-\varepsilon$ such that
\[
  |\{x\in E_{\varepsilon}:|Gf|\geq \lambda\}|\leq C_{\varepsilon}\biggr(\frac{\|f\|_p}{\lambda} \biggr)^q
\]
holds for all $\lambda >0$ and $f\in L^p (X)$, where $q=\min\{p,2\}$.
\end{Lemma}
When the operator $G$ in Lemma \ref{lem:s2t1} is a certain kind of maximal operator, its boundedness can be characterized by the following Lemma \ref{lem:s2t2}.
\begin{Lemma}[{\cite[Lemma, P.790]{Nikishin1972}}]\label{lem:s2t2}
Let $T_n:L^p (X)\rightarrow S[0,1]$ be a sequence of linear operators which are continuous in measure. If for each $f\in L^p (X)$ the $\lim_{n\rightarrow \infty}T_n f$ exists almost everywhere on $[0,1]$, then the operator $G$ defined by $Gf=\sup_n |T_n f|$ is hyperlinear and bounded.
\end{Lemma}
\begin{Remark}
Continuous in measure means that convergence of $f_k \rightarrow f_0$ in $L^p (X)$ implies convergence of $T_n f_k \rightarrow T_n f_0$ in measure on $[0,1]$. As for the Landau type Schr\"odinger operator $P_{a,\gamma}^t$, if $f_k \rightarrow f_0$ in $L^p (\mathbb R)$, where $1\leq p<2$, then by the Hausdorff-Young inequality, we get
\begin{align*}
&\quad \lim_{k\rightarrow \infty}|\{x\in [0,1]:|P_{a,\gamma}^t f_k(x)-P_{a,\gamma}^t f_0(x)|\geq \varepsilon\}|\\
&\leq \lim_{k\rightarrow \infty}\frac{1}{\varepsilon}\int_0^1 |P_{a,\gamma}^t f_k(x)-P_{a,\gamma}^t f_0(x)| \rmd x\\
&\lesssim \lim_{k\rightarrow \infty}\frac{1}{\varepsilon}\int_0^1 \biggr|\int_{\mathbb R}e^{i(x\xi+t|\xi|^a)}e^{-t^{\gamma}|\xi|^a}(\hat{f_k}-\hat{f_0})(\xi)\rmd \xi \biggr|\rmd x\\
&\leq \lim_{k\rightarrow \infty}\frac{1}{\varepsilon}\|e^{-t^{\gamma}|\cdot|^a}\|_p \cdot \|\hat{f_k}-\hat{f_0}\|_{p'}\lesssim \lim_{k\rightarrow \infty} \frac{1}{\varepsilon}\|f_k-f_0\|_p =0.
\end{align*}
Therefore, Theorem \ref{thm:s2t1} is an immediate consequence of Lemma \ref{lem:s2t1} and Lemma \ref{lem:s2t2}.
\end{Remark}

Next, inspired by the counterexample given in \cite[Theorem 1.3]{yuan2021dimension}, we give a proof of Theorem \ref{thm:s1t1}.\\
\textbf{Proof of Theorem \ref{thm:s1t1}.}
By contradiction, we assume that (\ref{eq:s1t2}) holds for all $f\in L^p (\mathbb R)$.
Then by Theorem \ref{thm:s2t1}, we obtain for any $\varepsilon >0$, there exists a set $E_{\varepsilon}\subset [0,1]$ with
$|E_{\varepsilon}|>1-\varepsilon$ such that the inequality (\ref{eq:s2t2}) holds for all $\lambda >0$ and $f\in L^p (\mathbb R)$.

The proof of Theorem \ref{thm:s1t1} is spilt in two cases:\\
\textbf{Case 1: $a\neq 1$.}
Let $g$ be a non-negative Schwartz function satisfying
\begin{equation}
g(\omega)=\begin{cases}
1,\quad |\omega|<\frac{1}{2};\\
0,\quad |\omega|>1.
\end{cases}
\nonumber
\end{equation}
For each $\theta\in (0,1/100)$ we define $\hat{f_{\theta}}(\omega)=\theta^{k_0}g(\theta^{k_0}\omega+1/\theta)$,
where $k_0 =2\lceil \gamma/(\min\{a,1\}(\gamma-1))\rceil$.
Then for any $x\in \mathbb R\setminus \{0\}$, by integration by parts, we have
\begin{align*}
f_{\theta}(x)
&=\frac{\theta^{k_0}}{2\pi}\int_{(-\theta-1)/\theta^{k_0+1}}^{(\theta-1)/\theta^{k_0+1}}
  g\biggr(\theta^{k_0}\omega+\frac{1}{\theta}\biggr)e^{ix\omega}\rmd\omega\\
&=-\frac{\theta^{2k_0}}{2\pi ix}\int_{(-\theta-1)/\theta^{k_0+1}}^{(\theta-1)/\theta^{k_0+1}}
  g'\biggr(\theta^{k_0}\omega+\frac{1}{\theta}\biggr)e^{ix\omega}\rmd\omega\\
&=-\frac{\theta^{3k_0}}{2\pi x^2}\int_{(-\theta-1)/\theta^{k_0+1}}^{(\theta-1)/\theta^{k_0+1}}
  g''\biggr(\theta^{k_0}\omega+\frac{1}{\theta}\biggr)e^{ix\omega}\rmd\omega.
\end{align*}
As a consequence of this, we obtain
\[
  |f_{\theta}(x)|\lesssim \frac{\theta^{3k_0}}{x^2}\int_{(-\theta-1)/\theta^{k_0+1}}^{(\theta-1)/\theta^{k_0+1}}
   \biggr|g''\biggr(\theta^{k_0}\omega+\frac{1}{\theta}\biggr)\biggr| \rmd\omega \lesssim \frac{\theta^{2k_0}}{x^2}.
\]
On the other hand, we have
\[
  |f_{\theta}(x)|\lesssim \theta^{k_0}\int_{(-\theta-1)/\theta^{k_0+1}}^{(\theta-1)/\theta^{k_0+1}}
   \biggr| g\biggr(\theta^{k_0}\omega+\frac{1}{\theta}\biggr)\biggr|\rmd\omega \lesssim 1.
\]
Therefore,
\begin{align}
\|f_{\theta}\|_p^p &=\int_{|x|\leq \theta^{k_0}}|f_{\theta}(x)|^p \rmd x+\int_{|x|> \theta^{k_0}}|f_{\theta}(x)|^p \rmd x \notag\\
&\lesssim \theta^{k_0} +\int_{\theta^{k_0}}^{+\infty} \frac{\theta^{2p k_0}}{x^{2p}}\rmd x\lesssim \theta^{k_0}.\label{eq:s2t3}
\end{align}
Recall that
\[
  |P_{a,\gamma}^t f_{\theta}(x)|\approx \biggr|\int_{\mathbb R}e^{i(x\omega+t|\omega|^a)}e^{-t^{\gamma}|\omega|^a}
   \theta^{k_0} g\biggr(\theta^{k_0}\omega+\frac{1}{\theta}\biggr)\rmd \omega \biggr|.
\]
Then if we make the change of variables $\xi=\theta^{k_0}\omega+1/\theta$, we get
\begin{align*}
|P_{a,\gamma}^t f_{\theta}(x)|
&\approx \biggr|\int_{\mathbb R}e^{i(x(\xi-1/\theta)/\theta^{k_0}+t|(\xi-1/\theta)/\theta^{k_0}|^a)}
  e^{-t^{\gamma}|(\xi-1/\theta)/\theta^{k_0}|^a}g(\xi)\rmd \xi \biggr|\\
&=\biggr|\int_{\mathbb R} e^{i(x\xi/\theta^{k_0}+t|(\xi-1/\theta)/\theta^{k_0}|^a)}
  e^{-t^{\gamma}|(\xi-1/\theta)/\theta^{k_0}|^a}g(\xi)\rmd \xi \biggr|.
\end{align*}
Let
\[
  \Phi_{x,t,\theta}(\xi)=x\frac{\xi}{\theta^{k_0}}+t\frac{|\theta\xi -1|^a}{\theta^{a(k_0 +1)}} -\frac{t}{\theta^{a(k_0 +1)}}
\]
and
\[
  \Psi_{t,\theta}(\xi)=t^{\gamma}\frac{|\theta\xi -1|^a}{\theta^{a(k_0 +1)}} .
\]
Note that $g$ is supported on $[-1,1]$, we obtain
\[
  |P_{a,\gamma}^t f_{\theta}(x)|\gtrsim \biggr|\int_{-1}^1 cos(\Phi_{x,t,\theta}(\xi)) e^{-\Psi_{t,\theta}(\xi)}g(\xi)\rmd \xi \biggr|.
\]
Next, we give estimates for $|\Phi_{x,t,\theta}(\xi)|$ and $|\Psi_{t,\theta}(\xi)|$ respectively. Applying Taylor's formula for $|\xi|\leq 1$, we have
\begin{equation}\label{eq:s2t4}
\frac{|\theta\xi -1|^a}{\theta^{a(k_0 +1)}}
=\frac{1}{\theta^{a(k_0 +1)}}-\frac{a\xi}{\theta^{a(k_0 +1)-1}}+\frac{a(a-1)}{2}\cdot \frac{\xi^2}{\theta^{a(k_0 +1)-2}}+R^{(2)}(\xi),
\end{equation}
where
\begin{equation}\label{eq:s2t5}
  R^{(2)}(\xi)=\frac{a(a-1)(a-2)(1-\theta \nu)^{a-3}(-\theta)^3 \xi^3}{6\cdot\theta^{a(k_0 +1)}}
\end{equation}
and $\nu$ is some number satisfying $|\nu|\leq |\xi|$. Then
\begin{equation}\label{eq:s2t6}
\Phi_{x,t,\theta}(\xi)=x\frac{\xi}{\theta^{k_0}}-\frac{ta\xi}{\theta^{a(k_0 +1)-1}}+\frac{a(a-1)}{2}\cdot \frac{t\xi^2}{\theta^{a(k_0 +1)-2}}+tR^{(2)}(\xi).
\end{equation}
Let $a_0 =\min\{a,1,1/|a-1|\}$. For any $x\in (0,a_0 \theta^{k_0 -1}]\subset [0,1]$, fix $t=x\theta^{a(k_0 +1)-1-k_0}/a$. It follows from $a(k_0 +1)-2> 2\gamma/(\gamma -1)-2>0$ that
\begin{equation}\label{eq:s2t7}
  0<t\leq \frac{a_0 \theta^{k_0 -1}\theta^{a(k_0 +1)-1-k_0}}{a}\leq \theta^{a(k_0 +1)-2}<1.
\end{equation}
Then by the equality (\ref{eq:s2t6}), we get
\begin{equation}\label{eq:s2t8}
\Phi_{x,t,\theta}(\xi)=\frac{a-1}{2}x\xi^2 \theta^{1-k_0}+tR^{(2)}(\xi).
\end{equation}
We claim that $|tR^{(2)}(\xi)|\leq 2/3$ holds for any $|\xi|\leq 1$. Indeed, according to (\ref{eq:s2t5}) and (\ref{eq:s2t7}),
\begin{align*}
|tR^{(2)}(\xi)| &\leq \theta^{a(k_0 +1)-2}\frac{|a(a-1)(a-2)|}{6\cdot\theta^{a(k_0 +1)}}(1-\theta \nu)^{a-3}\theta^3 \\
&=\frac{|a(a-1)(a-2)|}{6}(1-\theta \nu)^{a-3}\theta .
\end{align*}
It is easy to see that $(1-\theta \nu)^{a-3}<2^a (1-\theta \nu)^{-3}\leq 2^{a+1}$ due to $|\nu|\leq 1$ and $\theta<1/100$.
Therefore, when $\theta$ is sufficiently small, we have
\begin{equation}\label{eq:s2t9}
|tR^{(2)}(\xi)|\leq \frac{|a(a-1)(a-2)|}{3}2^a \theta <\frac{2}{3}.
\end{equation}
Combining (\ref{eq:s2t8}) and (\ref{eq:s2t9}), we obtain
\begin{equation}\label{eq:s2t10}
|\Phi_{x,t,\theta}(\xi)|\leq \frac{|a-1|}{2}\cdot\frac{1}{|a-1|}\theta^{k_0 -1}\theta^{1-k_0}+\frac{2}{3}<\frac{\pi}{2}.
\end{equation}
Similarly, by equality (\ref{eq:s2t4}),
\[
  \Psi_{t,\theta}(\xi)=\frac{t^{\gamma}}{\theta^{a(k_0 +1)}}-\frac{t^{\gamma}a\xi}{\theta^{a(k_0 +1)-1}}
  +\frac{a(a-1)}{2}\cdot \frac{t^{\gamma}\xi^2}{\theta^{a(k_0 +1)-2}}+t^{\gamma}R^{(2)}(\xi).
\]
Note that $|t^{\gamma}R^{(2)}(\xi)|\leq|tR^{(2)}(\xi)|\leq 2/3$ and $|\xi|\leq \theta^{-1}$, we have
\begin{align*}
|\Psi_{t,\theta}(\xi)|
&\lesssim \frac{t^{\gamma}}{\theta^{a(k_0 +1)}}+\frac{t^{\gamma}}{\theta^{a(k_0 +1)}}\cdot \frac{|\xi|}{\theta^{-1}}
 +\frac{t^{\gamma}}{\theta^{a(k_0 +1)}}\cdot \frac{\xi^2}{\theta^{-2}}+|t^{\gamma}R^{(2)}(\xi)|\\
&\lesssim \frac{t^{\gamma}}{\theta^{a(k_0 +1)}}+1\\
&\lesssim x^{\gamma}\theta^{a\gamma(k_0 +1)-\gamma-\gamma k_0}\theta^{-a(k_0 +1)}+1\\
&\lesssim \theta^{(k_0 -1)\gamma}\theta^{a\gamma(k_0 +1)-\gamma-\gamma k_0}\theta^{-a(k_0 +1)}+1\\
&= \theta^{a(k_0 +1)(\gamma-1)-2\gamma}+1.
\end{align*}
Since $a(k_0 +1)(\gamma-1)-2\gamma>2a(\gamma-1)\gamma/a(\gamma-1)-2\gamma=0$, we have
\begin{equation}\label{eq:s2t11}
|\Psi_{t,\theta}(\xi)|\lesssim \theta^{a(k_0 +1)(\gamma-1)-2\gamma}+1\lesssim 1.
\end{equation}
Then, combining (\ref{eq:s2t10}) and (\ref{eq:s2t11}), we get $P_{a,\gamma}^* f_{\theta} (x)\gtrsim 1$ holds for all $x\in (0,a_0 \theta^{k_0 -1}]$. If we take $\varepsilon =1/4$ in (\ref{eq:s2t2}), then we can see from the proof of {\cite[Theorem 1.3]{yuan2021dimension}} that for any $\lambda\in (0,1/10)$, there exists $x_0\in [0,1]$ such that
\begin{equation}\label{eq:s2t12}
  |E_{\varepsilon}\cap [x_0 ,x_0+\lambda]|\geq \frac{1}{2}\lambda.
\end{equation}
Let $\lambda =a_0 \theta^{k_0 -1}$ and $\widetilde{f_{\theta}}(x)=f_{\theta}(x-x_0)$. Therefore, by inequalities (\ref{eq:s2t2}), (\ref{eq:s2t3}) and (\ref{eq:s2t12}),
\begin{align*}
\frac{1}{2}a_0 \theta^{k_0 -1}&\leq |E_{\varepsilon}\cap [x_0 ,x_0+a_0 \theta^{k_0 -1}]|\\
&=|\{x\in E_{\varepsilon}\cap [x_0 ,x_0+a_0 \theta^{k_0 -1}]:P_{a,\gamma}^* \widetilde{f_{\theta}} (x)\gtrsim 1\}|\\
&\leq |\{x\in E_{\varepsilon}:P_{a,\gamma}^* \widetilde{f_{\theta}} (x)\gtrsim 1\}|\\
&\lesssim \|\widetilde{f_{\theta}}\|_p^p =\|f_{\theta}\|_p^p \lesssim \theta^{k_0}.
\end{align*}
That is, $\theta\gtrsim 1$. Letting $\theta\rightarrow 0$ gives a contradiction.\\
\textbf{Case 2: $a= 1$.}
Let $\theta\in (0,1/100)$ and $g$ be defined as in the previous case.
Define $\hat{f_{\theta}}(\omega)=\theta^{k_0}g(\theta^{k_0}\omega+1/\theta)$, where $k_0 =2\lceil \gamma/(\gamma-1)\rceil$.
Then from the proof of case $1$, we see that $\|f_\theta\|_p^p \lesssim \theta^{k_0}$.
Since $\xi< \theta^{-1}$, we have
\begin{align*}
|P_{a,\gamma}^t f_{\theta}(x)|
&\approx \biggr|\int_{\mathbb R}e^{i(x(\xi-1/\theta)/\theta^{k_0}+t|(\xi-1/\theta)/\theta^{k_0}|)}
  e^{-t^{\gamma}|(\xi-1/\theta)/\theta^{k_0}|}g(\xi)\rmd \xi \biggr|\\
&=\biggr|\int_{\mathbb R} e^{i(x\xi/\theta^{k_0}-t\xi/\theta^{k_0})}
  e^{-t^{\gamma} (1/\theta^{k_0 +1}-\xi/\theta^{k_0})}g(\xi)\rmd \xi \biggr|.
\end{align*}
Let
\[
  \Phi_{x,t,\theta}(\xi)=x\frac{\xi}{\theta^{k_0}}-t\frac{\xi}{\theta^{k_0}}
\]
and
\[
  \Psi_{t,\theta}(\xi)=t^{\gamma}\biggr ( \frac{1}{\theta^{k_0 +1}}-\frac{\xi}{\theta^{k_0}} \biggr).
\]
Therefore, we have
\[
  |P_{a,\gamma}^t f_{\theta}(x)|\gtrsim \biggr|\int_{-1}^1 cos(\Phi_{x,t,\theta}(\xi)) e^{-\Psi_{t,\theta}(\xi)}g(\xi)\rmd \xi \biggr|.
\]
For any $x\in (0,\theta^{k_0 -1}]$, fix $t=x\subset [0,1]$. Then we have $\Phi_{x,t,\theta}(\xi)=0$. Besides,
\begin{align*}
|\Psi_{t,\theta}(\xi)|&\leq \frac{t^{\gamma}}{\theta^{k_0 +1}}+\frac{t^{\gamma}|\xi|}{\theta^{k_0}}\lesssim \frac{t^{\gamma}}{\theta^{k_0 +1}}\\
&\leq \theta^{(k_0 -1)\gamma-k_0 -1}=\theta^{(\gamma -1)k_0 -\gamma -1}.
\end{align*}
Note that $(\gamma -1)k_0 -\gamma -1\geq \gamma -1>0$, we obtain $|\Psi_{t,\theta}(\xi)|\lesssim 1$.
Thus $P_{a,\gamma}^* f_{\theta} (x)\gtrsim 1$ holds for all $x\in (0,\theta^{k_0 -1}]$.
Similar to the previous argument, we have
\begin{align*}
\frac{1}{2}\theta^{k_0 -1}&\leq |E_{\varepsilon}\cap [x_0 ,x_0+\theta^{k_0 -1}]|\\
&\leq |\{x\in E_{\varepsilon}:P_{a,\gamma}^* \widetilde{f_{\theta}} (x)\gtrsim 1\}|\\
&\lesssim \|\widetilde{f_{\theta}}\|_p^p =\|f_{\theta}\|_p^p \lesssim \theta^{k_0},
\end{align*}
which gives the desired contradiction if we take $\theta\rightarrow 0$.\hfill $\Box$

\subsection{Proof of Theorem \ref{thm:s1t2}}

To obtain the pointwise convergence result of the Landau type Schr\"odinger operator, we first prove the following Lemma \ref{lem:s3t1}.
\begin{Lemma}\label{lem:s3t1}
Let $s>0$ and $1<p\leq 2$. Then for any $u\in W^{s,p}(\mathbb R)$, we have
\begin{equation}\label{eq:s3t1}
\|(1+|\cdot|^2)^{s/2}\hat{u}\|_{L^{p'}(\mathbb R)}\lesssim \|u\|_{W^{s,p}(\mathbb R)},
\end{equation}
where $p'$ is the conjugate exponent of $p$.
\end{Lemma}

In fact, for the case $0<s<1$, Lemma \ref{lem:s3t1} is an immediate consequence of \cite[Theorem 5(c), P.155]{stein1970singular},
whose proof depends on the Littlewood-Paley theorem. In this subsection, we provide an alternative proof.
\proof We first consider the case $0<s<1$. We claim that
\begin{equation}\label{eq:s3t2}
  \biggr(\int_{\mathbb R}|\hat{u}(\xi)|^{p'} |\xi|^{sp'}\rmd \xi\biggr)^{1/p'}
\lesssim \biggr(\iint_{\mathbb R\times \mathbb R}\frac{|u(x)-u(y)|^p}{|x-y|^{sp+1}}\rmd x\rmd y\biggr)^{1/p}.
\end{equation}
Indeed, by the Hausdorff-Young inequality,
\begin{align*}
\iint_{\mathbb R \times \mathbb R}\frac{|u(x)-u(y)|^p}{|x-y|^{sp+1}}\rmd x\rmd y
&=\int_{\mathbb R}\int_{\mathbb R}\frac{|u(x+y)-u(y)|^p}{|x|^{sp+1}}\rmd y\rmd x\\
&\geq \biggr(\int_{\mathbb R}\biggr(\int_{\mathbb R}\frac{|(u-u(\cdot+x))\, \hat{}(y)|^{p'}}{|x|^{sp'+p'/p}}\rmd y \biggr)^{p/p'}\rmd x\biggr)^{p/p}.
\end{align*}
It is obvious that $p'\geq p$ due to $1<p\leq 2$. Then by the Minkowski's integral inequality, we get
\begin{align*}
\iint_{\mathbb R \times \mathbb R}\frac{|u(x)-u(y)|^p}{|x-y|^{sp+1}}\rmd x\rmd y
&\geq \biggr(\int_{\mathbb R}\biggr(\int_{\mathbb R}\frac{|(u-u(\cdot+x))\, \hat{}(y)|^p}{|x|^{sp+1}}\rmd x \biggr)^{p'/p}\rmd y\biggr)^{p/p'}\\
&= \biggr(\int_{\mathbb R}\biggr(\int_{\mathbb R}\frac{|\hat{u}(y)(1-e^{ixy})|^p}{|x|^{sp+1}}\rmd x \biggr)^{p'/p}\rmd y\biggr)^{p/p'}\\
&= \biggr(\int_{\mathbb R}\biggr(\int_{\mathbb R}\frac{|2\sin (xy/2)|^p}{|x|^{sp+1}}\rmd x \biggr)^{p'/p}|\hat{u}(y)|^{p'}\rmd y\biggr)^{p/p'}\\
&\gtrsim \biggr(\int_{\mathbb R}|\hat{u}(y)|^{p'}|y|^{sp'}\rmd y \biggr)^{p/p'}.
\end{align*}
Hence, we arrive at (\ref{eq:s3t2}). Applying the Hausdorff-Young inequality again, we have
\begin{align*}
&\quad \int_{\mathbb R}|u(y)|^p \rmd y+\iint_{\mathbb R \times \mathbb R}\frac{|u(x)-u(y)|^p}{|x-y|^{sp+1}}\rmd x\rmd y\\
&\gtrsim \biggr(\int_{\mathbb R}|\hat{u}(y)|^{p'}\rmd y \biggr)^{p/p'}+\biggr(\int_{\mathbb R}|\hat{u}(y)|^{p'}|y|^{sp'}\rmd y \biggr)^{p/p'}\\
&\geq \biggr(\int_{\mathbb R}|\hat{u}(y)|^{p'}(1+|y|^{sp'})\rmd y \biggr)^{p/p'}\\
&\gtrsim \biggr(\int_{\mathbb R}(1+|y|^2)^{sp'/2}|\hat{u}(y)|^{p'}\rmd y \biggr)^{p/p'}.
\end{align*}
That is, $\|(1+|\cdot|^2)^{s/2}\hat{u}\|_{L^{p'}(\mathbb R)}\lesssim \|u\|_{W^{s,p}(\mathbb R)}$.

Next, we consider the case where $s=k\geq 1$ be an integer. Note that the distributional derivative $D^k u$ satisfies $|(D^k u)\, \hat{}(\xi)|=|\xi^k \hat{u}(\xi)|$. Therefore, we obtain
\begin{align*}
\|(1+|\cdot|^2)^{k/2}\hat{u}\|_{L^{p'}}
&=\biggr(\int_{\mathbb R}(1+|\xi|^2)^{kp'/2}|\hat{u}(\xi)|^{p'}\rmd \xi \biggr)^{1/p'}\\
&\lesssim \biggr(\int_{|\xi|\leq 1}|\hat{u}(\xi)|^{p'}\rmd \xi+ \int_{|\xi|> 1}|\xi|^{kp'}|\hat{u}(\xi)|^{p'}\rmd \xi\biggr)^{1/p'}\\
&\leq \biggr(\int_{\mathbb R}|\hat{u}(\xi)|^{p'}\rmd \xi+ \int_{\mathbb R}|(D^k u)\,\hat{}(\xi)|^{p'}\rmd \xi\biggr)^{1/p'}.
\end{align*}
Then by the Hausdorff-Young inequality, we have
\begin{align*}
\|(1+|\cdot|^2)^{k/2}\hat{u}\|_{L^{p'}}
&\lesssim \biggr(\biggr(\int_{\mathbb R}|u(x)|^p \rmd x\biggr)^{p'/p}+\biggr(\int_{\mathbb R}|(D^k u)(x)|^p \rmd x\biggr)^{p'/p} \biggr)^{1/p'}\\
&\leq (\|u\|_{L^p (\mathbb R)}^p +\|D^k u\|_{L^p (\mathbb R)}^p)^{1/p}\leq \|u\|_{W^{k,p}(\mathbb R)}.
\end{align*}

Finally, we consider the case $s>1$ not be an integer.
For this case, we decompose $s$ as $s=m+\tau$, where $m\geq 1$ be an integer and $\tau\in (0,1)$.
It is easy to see that
\[
  (1+|\xi|^2)^{(m+\tau)p'/2}\leq 2^{\tau p'/2}(1+|\xi|^2)^{mp'/2}+2^{mp'/2}|\xi|^{mp'}(1+|\xi|^2)^{\tau p'/2}
\]
and $|(D^m u)\, \hat{}(\xi)|=|\xi^m \hat{u}(\xi)|$. Hence, combining the two cases we have proved, we get
\begin{align*}
&\qquad \|(1+|\cdot|^2)^{s/2}\hat{u}\|_{L^{p'}}\\
&=\biggr(\int_{\mathbb R}(1+|\xi|^2)^{(m+\tau)p'/2}|\hat{u}(\xi)|^{p'}\rmd \xi \biggr)^{1/p'}\\
&\lesssim \biggr(\int_{\mathbb R}(1+|\xi|^2)^{mp'/2}|\hat{u}(\xi)|^{p'}\rmd \xi
+\int_{\mathbb R}(1+|\xi|^2)^{\tau p'/2}|(D^m u)\,\hat{}(\xi)|^{p'}\rmd \xi\biggr)^{1/p'}\\
&\lesssim \biggr(\|u\|_{W^{m,p}(\mathbb R)}^{p'}
+\biggr(\|D^m u\|_{L^p (\mathbb R)}^p +\iint_{\mathbb R \times \mathbb R}\frac{|D^m u(x)-D^m u(y)|^p}{|x-y|^{\tau p+1}}\rmd x\rmd y \biggr)^{p'/p} \biggr)^{1/p'}\\
&\lesssim \biggr(\|u\|_{W^{m,p}(\mathbb R)}^p +\iint_{\mathbb R \times \mathbb R}\frac{|D^m u(x)-D^m u(y)|^p}{|x-y|^{\tau p+1}}\rmd x\rmd y\biggr)^{1/p}
=\|u\|_{W^{s,p}(\mathbb R)},
\end{align*}
which completes the proof of Lemma \ref{lem:s3t1}.
\endproof

\begin{Remark}
In fact, for the case of $p=2$ and $0<s<1$, we can further prove that the two sides of (\ref{eq:s3t1}) are equivalent.
For a detailed proof we refer to {\cite[Proposition 3.4]{Nezza2012}}.
\end{Remark}

Now we are ready to prove Theorem \ref{thm:s1t2}.\\
\textbf{Proof of Theorem \ref{thm:s1t2}.}
Take some function $\psi$ such that $\hat{\psi}\in C_c^\infty (\mathbb R)$, $0\leq \hat{\psi}(\omega)\leq 1$, and
\begin{equation}
\hat{\psi}(\omega)=\begin{cases}
1,\quad |\omega|<1;\\
0,\quad |\omega|\geq 2.
\end{cases}
\nonumber
\end{equation}
Set $f_1 =f * \psi$, $f_2 =f-f_1$. Then by the Hausdorff-Young inequality, we obtain
\[
  \|\hat{f_1}\|_{L^1}=\|\hat{f}\cdot \hat{\psi}\|_{L^1}\leq \|\hat{f}\|_{L^{p'}}\|\hat{\psi}\|_{L^p}\leq  \|f\|_{L^p}\|\hat{\psi}\|_{L^p}<\infty .
\]
Applying dominated convergence theorem yields
\begin{equation}\label{eq:s3t3}
\lim_{t\rightarrow 0}P_{a,\gamma}^t f_1 (x)= f_1 (x),\quad a.e. \ x\in \mathbb R.
\end{equation}
Now we first consider the case $a>1$. Let $1/r= 1/p-1/2$. Then we have $1/r<s-\min\{a(1-1/\gamma)/4,1/4\}$.
Hence, there is some $\sigma >\min\{a(1-1/\gamma)/4,1/4\}$ such that $1/r<s-\sigma$.
Consequently, by H\"older's inequality,
\begin{align}
\|f_2\|_{H^{\sigma}}
&= \biggr(\int_{|\omega|>1}(1+|\omega|^2)^{\sigma}(1-\hat{\psi}(\omega))^2 |\hat{f}(\omega)|^2 \rmd \omega \biggr)^{1/2}
\lesssim \||\omega|^{\sigma}\hat{f}(\omega)\chi_{\{|\omega|>1\}}\|_{L^2} \notag\\
&\leq \||\omega|^{-(s-\sigma)}\chi_{\{|\omega|>1\}}\|_{L^r} \||\omega|^s \hat{f}(\omega)\|_{L^{p'}}
\lesssim \||\omega|^s \hat{f}(\omega)\|_{L^{p'}}. \label{eq:s3t4}
\end{align}
It follows from Lemma \ref{lem:s3t1} that
\[
  \||\omega|^s \hat{f}(\omega)\|_{L^{p'}}\leq \|(1+|\omega|^2)^{s/2} \hat{f}(\omega)\|_{L^{p'}}\lesssim \|f\|_{W^{s,p}(\mathbb R)}<\infty.
\]
Thus $f_2 \in H^{\sigma}$.
Recall that it was proved in {\cite[Theorem 1.3]{bailey2013boundedness}} that $\lim_{t\rightarrow 0}P_{a,\gamma}^t f (x)= f (x),\ a.e. \ x\in \mathbb R$
for all $f\in H^{\sigma}$ with $\sigma >\min\{a(1-1/\gamma)/4,1/4\}$. Therefore,
\begin{equation}\label{eq:s3t5}
\lim_{t\rightarrow 0}P_{a,\gamma}^t f_2 (x)= f_2 (x),\quad a.e. \ x\in \mathbb R.
\end{equation}
Hence, (\ref{eq:s1t2}) follows immediately by (\ref{eq:s3t3}) and (\ref{eq:s3t5}).
For the case $0<a\leq 1$, it can be shown by {\cite[Corollary 1.2]{yuan2021dimension}} that
(\ref{eq:s1t2}) still holds under the assumptions of Theorem \ref{thm:s1t2}.
The proof is almost the same and so is omitted. \hfill $\Box$

\subsection{Proof of Theorem \ref{thm:s1t3}.}
The following lemma will be crucial for the maximal estimate in the proof of Theorem \ref{thm:s1t3}.
\begin{Lemma}[{\cite[Lemma 5.1]{Li2021JFA}}]\label{lem:s4t1}
Let $P(\xi)$ be a real continuous function. Assume that $g$ is a Schwartz function whose Fourier transform is supported in the annulus
$A(\lambda)=\{\xi\in \mathbb R^n :|\xi|\approx \lambda\}$. Suppose $\alpha >0$ and $\gamma (x,t)$ satisfy that
\[
  |\gamma (x,t)-x|\lesssim t^{\alpha},\quad \gamma(x,0)=x
\]
uniformly for all $x\in B(x_0,R)$ and $t\in (0,\lambda^{-1/\alpha})$. Then for each $x\in B(x_0,R)$ and $t\in (0,\lambda^{-1/\alpha})$, we have
\begin{align*}
|e^{itP(D)}g(\gamma(x,t))|
&:=\biggr|\int_{\mathbb R^n}e^{i\gamma(x,t)\cdot \xi}e^{itP(\xi)}\hat{g}(\xi)\rmd \xi \biggr|\\
&\leq \sum_{l\in \mathbb Z^n}\frac{C_n}{(1+|l|)^{n+1}}\biggr|\int_{\mathbb R^n}e^{i(x+l/\lambda)\cdot \xi+itP(\xi)}\hat{g}(\xi)\rmd \xi\biggr|.
\end{align*}
\end{Lemma}

Next, inspired by the work of Li and Wang \cite{Li2021arxiv}, we give a proof of Theorem \ref{thm:s1t3}.\\
\textbf{Proof of Theorem \ref{thm:s1t3}.}
$(1)$ We claim that it suffices to show that for some $q\geq 1$ and any $\varepsilon >0$, $s_1=s_0+\varepsilon$, we have
\begin{equation}\label{eq:s4t1}
\biggr\|\sup_{0<t<1}\frac{|P_{a,\gamma}^t (f)(\Gamma(x,t))-f(x)|}{t^{\delta\min\{1,\gamma\}/a}}\biggr\|_{L^q (B(x_0 ,R))}
\leq C_{\varepsilon}\|f\|_{H^{s_1+\delta}(\mathbb R)},
\end{equation}
where $0<\delta<a\min\{\beta,\gamma\}/\min\{1,\gamma\}$. Assume at the moment that (\ref{eq:s4t1}) holds. Then given $\lambda >0$,
there exists some $g\in C_c^\infty(\mathbb R)$ such that
\[
  \|f-g\|_{H^{s_1+\delta}(\mathbb R)}\leq \frac{\lambda \varepsilon^{1/q}}{2C_{\varepsilon}},
\]
which follows
\begin{align}
&\quad \biggr|\biggr\{x\in B(x_0,R):\sup_{0<t<1}\frac{|P_{a,\gamma}^t (f-g)(\Gamma(x,t))-(f-g)(x)|}{t^{\delta\min\{1,\gamma\}/a}}
>\frac{\lambda}{2}\biggr\} \biggr| \notag\\
&\leq \frac{2^q}{\lambda^q}\biggr\|\sup_{0<t<1}\frac{|P_{a,\gamma}^t (f-g)(\Gamma(x,t))-(f-g)(x)|}{t^{\delta\min\{1,\gamma\}/a}}
       \biggr\|_{L^q (B(x_0,R))}^q \notag\\
&\leq \frac{2^q C_{\varepsilon}^q}{\lambda^q} \|f-g\|_{H^{s_1+\delta}(\mathbb R)}^q \leq \varepsilon.\label{eq:s4t2}
\end{align}
Furthermore, we have
\begin{equation}\label{eq:s4t3}
\frac{|P_{a,\gamma}^t (g)(\Gamma(x,t))-g(x)|}{t^{\delta\min\{1,\gamma\}/a}}\rightarrow 0,\quad if \quad t\rightarrow 0^+
\end{equation}
uniformly for $x\in B(x_0,R)$. Indeed, for each $x\in B(x_0,R)$,
\begin{align}
\lim_{t\rightarrow 0^+}\frac{|P_{a,\gamma}^t (g)(\Gamma(x,t))-g(x)|}{t^{\delta\min\{1,\gamma\}/a}}
&\leq \lim_{t\rightarrow 0^+}\frac{|P_{a,\gamma}^t (g)(\Gamma(x,t))-g(\Gamma(x,t))|}{t^{\delta\min\{1,\gamma\}/a}}\notag\\
&+ \lim_{t\rightarrow 0^+}\frac{|g(\Gamma(x,t))-g(x)|}{t^{\delta\min\{1,\gamma\}/a}}.\label{eq:s4t4}
\end{align}
It is obvious that $|e^{-y}-1|\leq y$ holds for all $y\geq 0$. Therefore,
\begin{align*}
|e^{it|\xi|^a}e^{-t^{\gamma}|\xi|^a}-1|
&= |e^{it|\xi|^a}e^{-t^{\gamma}|\xi|^a}-e^{-t^{\gamma}|\xi|^a}+e^{-t^{\gamma}|\xi|^a}-1|\\
&\leq |e^{it|\xi|^a}-1|+|e^{-t^{\gamma}|\xi|^a}-1|\\
&\leq t|\xi|^a +t^{\gamma}|\xi|^a \lesssim t^{\min\{1,\gamma\}}|\xi|^a
\end{align*}
holds for all $t\in [0,1]$. As a consequence of this, we get
\begin{align*}
|P_{a,\gamma}^t (g)(\Gamma(x,t))-g(\Gamma(x,t))|
&\approx \biggr| \int_{\mathbb R}e^{i\Gamma(x,t)\cdot \xi}(e^{it|\xi|^a}e^{-t^{\gamma}|\xi|^a}-1)\hat{g}(\xi)\rmd \xi\biggr|\\
&\lesssim t^{\min\{1,\gamma\}}\int_{\mathbb R}|\xi|^a |\hat{g}(\xi)|\rmd \xi.
\end{align*}
Consequently,
\begin{equation}\label{eq:s4t5}
\frac{|P_{a,\gamma}^t (g)(\Gamma(x,t))-g(\Gamma(x,t))|}{t^{\delta\min\{1,\gamma\}/a}}
\lesssim t^{\min\{1,\gamma\}(1-\delta/a)}\int_{\mathbb R}|\xi|^a |\hat{g}(\xi)|\rmd \xi.
\end{equation}
On the other hand,
\begin{align}
\frac{|g(\Gamma(x,t))-g(x)|}{t^{\delta\min\{1,\gamma\}/a}}
&\lesssim t^{-\delta\min\{1,\gamma\}/a}\int_{\mathbb R}|e^{i(\Gamma(x,t)-x)\cdot\xi}-1|\cdot|\hat{g}(\xi)|\rmd \xi \notag\\
&\leq \frac{|\Gamma(x,t)-x|}{t^{\delta\min\{1,\gamma\}/a}}\int_{\mathbb R}|\xi|\cdot|\hat{g}(\xi)|\rmd \xi \notag\\
&\leq t^{\beta-\delta\min\{1,\gamma\}/a}\int_{\mathbb R}|\xi|\cdot|\hat{g}(\xi)|\rmd \xi.\label{eq:s4t6}
\end{align}
Note that $\delta <a\min\{\beta,\gamma\}/\min\{1,\gamma\}$. Hence, inequalities (\ref{eq:s4t4})-(\ref{eq:s4t6}) imply (\ref{eq:s4t3}).
By (\ref{eq:s4t2}) and (\ref{eq:s4t3}), we obtain
\[
  \biggr|\biggr\{x\in B(x_0,R):\limsup_{t\rightarrow 0^+}
   \frac{|P_{a,\gamma}^t f(\Gamma(x,t))-f(x)|}{t^{\delta\min\{1,\gamma\}/a}}>\lambda\biggr\} \biggr| \leq \varepsilon,
\]
which indicates that (\ref{eq:s1t5}) holds for all $f\in H^{s_1+\delta}(\mathbb R)$ and almost every $x\in B(x_0,R)$.
By the arbitrariness of $\varepsilon$, we obtain (\ref{eq:s1t5}) for all $f\in H^{s+\delta}(\mathbb R)$, $s>s_0$.

Next, we will prove (\ref{eq:s4t1}) for $q=\min\{p,2\}$. In order to prove (\ref{eq:s4t1}), we decompose $f$ as follows:
\[
  f=\sum_{k=0}^\infty f_k,
\]
where $supp\hat{f_0}\subset B(0,1)$, $supp \hat{f_k}\subset \{\xi:|\xi|\approx 2^k\}$, $k\geq 1$. It follows that
\begin{align}
&\quad \biggr\|\sup_{0<t<1}\frac{|P_{a,\gamma}^t (f)(\Gamma(x,t))-f(x)|}{t^{\delta\min\{1,\gamma\}/a}}\biggr\|_{L^q (B(x_0 ,R))} \notag\\
&\leq \sum_{k=0}^\infty \biggr\|\sup_{0<t<1}\frac{|P_{a,\gamma}^t (f_k)(\Gamma(x,t))-f_k(x)|}{t^{\delta\min\{1,\gamma\}/a}}\biggr\|_{L^q (B(x_0 ,R))}.\label{eq:s4t7}
\end{align}
For the case $k\lesssim 1$, there exists a constant $M>0$ such that $supp \hat{f_k}\subset \{\xi:|\xi|\leq M\}$.
By invoking analogous arguments to that in (\ref{eq:s4t5}) and (\ref{eq:s4t6}), we have
\begin{align}
&\quad \biggr\|\sup_{0<t<1}\frac{|P_{a,\gamma}^t (f_k)(\Gamma(x,t))-f_k(x)|}{t^{\delta\min\{1,\gamma\}/a}}\biggr\|_{L^q (B(x_0 ,R))} \notag\\
&\leq \biggr\|\sup_{0<t<1}\frac{|P_{a,\gamma}^t (f_k)(\Gamma(x,t))-f_k(\Gamma(x,t))|}{t^{\delta\min\{1,\gamma\}/a}}\biggr\|_{L^q (B(x_0 ,R))}
 +\biggr\|\sup_{0<t<1}\frac{|f_k(\Gamma(x,t))-f_k(x)|}{t^{\delta\min\{1,\gamma\}/a}}\biggr\|_{L^q (B(x_0 ,R))} \notag\\
&\lesssim \|f\|_{H^{s_1+\delta}(\mathbb R)}.\label{eq:s4t8}
\end{align}
For the case $k\gg 1$,
\begin{align}
&\quad \biggr\|\sup_{0<t<1}\frac{|P_{a,\gamma}^t (f_k)(\Gamma(x,t))-f_k(x)|}{t^{\delta\min\{1,\gamma\}/a}}\biggr\|_{L^q (B(x_0 ,R))} \notag\\
&\leq \biggr\|\sup_{2^{-ak/\min\{1,\gamma\}}\leq t<1}
  \frac{|P_{a,\gamma}^t (f_k)(\Gamma(x,t))-f_k(x)|}{t^{\delta\min\{1,\gamma\}/a}}\biggr\|_{L^q (B(x_0 ,R))} \notag\\
&\qquad +\biggr\|\sup_{0<t<2^{-ak/\min\{1,\gamma\}}}
  \frac{|P_{a,\gamma}^t (f_k)(\Gamma(x,t))-f_k(x)|}{t^{\delta\min\{1,\gamma\}/a}}\biggr\|_{L^q (B(x_0 ,R))} :=I+II.\label{eq:s4t9}
\end{align}
We first estimate the part $I$. From (\ref{eq:s1t4}) we obtain
\[
  \biggr\|\sup_{2^{-ak/\min\{1,\gamma\}}\leq t<1}|P_{a,\gamma}^t (f_k)(\Gamma(x,t))|\biggr\|_{L^p (B(x_0 ,R))}
   \lesssim 2^{(s_0+\varepsilon/2)k}\|f_k\|_{L^2 (\mathbb R)}.
\]
Recall that $q=\min\{p,2\}$, we have
\begin{align}
I &\leq 2^{\delta k}\biggr\|\sup_{2^{-ak/\min\{1,\gamma\}}\leq t<1}|P_{a,\gamma}^t (f_k)(\Gamma(x,t))-f_k(x)|\biggr\|_{L^q (B(x_0 ,R))} \notag\\
&\lesssim 2^{\delta k} \biggr\{\biggr\| \sup_{2^{-ak/\min\{1,\gamma\}}\leq t<1}|P_{a,\gamma}^t (f_k)(\Gamma(x,t))| \biggr\|_{L^p (B(x_0 ,R))}
 +\|f_k\|_{L^2 (B(x_0,R))}\biggr\} \notag\\
&\lesssim 2^{(\delta +s_0+\varepsilon/2)k}\|f_k\|_{L^2 (\mathbb R)} \notag\\
&\lesssim 2^{-\varepsilon k/2}\|f\|_{H^{s_1+\delta}(\mathbb R)}.\label{eq:s4t10}
\end{align}
For the part $II$, by Minkowski's inequality,
\begin{align}
II &\leq \biggr\|\sup_{0<t<2^{-ak/\min\{1,\gamma\}}}
     \frac{|P_{a,\gamma}^t (f_k)(\Gamma(x,t))-f_k(\Gamma(x,t))|}{t^{\delta\min\{1,\gamma\}/a}}\biggr\|_{L^q (B(x_0 ,R))} \notag\\
&\qquad +\biggr\|\sup_{0<t<2^{-ak/\min\{1,\gamma\}}}
  \frac{|f_k(\Gamma(x,t))-f_k(x)|}{t^{\delta\min\{1,\gamma\}/a}}\biggr\|_{L^q (B(x_0 ,R))}.\label{eq:s4t11}
\end{align}
By Taylor's formula, we get
\begin{align}
&\quad \biggr\|\sup_{0<t<2^{-ak/\min\{1,\gamma\}}}\frac{|f_k(\Gamma(x,t))-f_k(x)|}{t^{\delta\min\{1,\gamma\}/a}}\biggr\|_{L^q (B(x_0 ,R))}\notag\\
&\lesssim \sum_{j\geq 1}\frac{1}{j!}\biggr\|\sup_{0<t<2^{-ak/\min\{1,\gamma\}}}\frac{|\Gamma(x,t)-x|^j}{t^{\delta\min\{1,\gamma\}/a}}
  \biggr| \int_{\mathbb R}e^{ix\cdot \xi}\xi^j \hat{f_k}(\xi) \rmd\xi\biggr| \biggr\|_{L^2 (B(x_0,R))} \notag\\
&\lesssim \sum_{j\geq 1}\frac{1}{j!}2^{-ak\beta j/\min\{1,\gamma\}}2^{\delta k}
  \biggr\|\int_{\mathbb R}e^{ix\cdot \xi} \xi^j \hat{f_k}(\xi) \rmd\xi \biggr\|_{L^2 (B(x_0,R))} \notag\\
&\lesssim \sum_{j\geq 1}\frac{1}{j!}2^{-ak\beta j/\min\{1,\gamma\}+\delta k}2^{kj} \|\hat{f_k}\|_{L^2 (\mathbb R)} \notag\\
&\lesssim 2^{\delta k}\|\hat{f_k}\|_{L^2 (\mathbb R)}\lesssim 2^{-s_1 k}\|f\|_{H^{s_1+\delta}(\mathbb R)}.\label{eq:s4t12}
\end{align}
Applying Taylor's formula again, we obtain
\begin{align}
|P_{a,\gamma}^t (f_k)(\Gamma(x,t))-f_k(\Gamma(x,t))|
&\leq \sum_{j=1}^\infty \frac{1}{j!}\biggr|\int_{\mathbb R}e^{i\Gamma(x,t)\cdot \xi}(it|\xi|^a -t^{\gamma}|\xi|^a)^j \hat{f_k}(\xi)\rmd \xi\biggr| \notag\\
&\lesssim \sum_{j=1}^\infty \frac{2^{j/2}}{j!}t^{\min\{1,\gamma\}j}
  \biggr|\int_{\mathbb R}e^{i\Gamma(x,t)\cdot \xi}|\xi|^{aj}\hat{f_k}(\xi)\rmd \xi\biggr|.\label{eq:s4t13}
\end{align}
It follows from $\beta\geq \min\{1,\gamma\}/a$ that $2^{-ak/\min\{1,\gamma\}}\leq 2^{-k/\beta}$.
If we take $\lambda =2^k$, $\alpha=\beta$ and $P(\xi)\equiv 0$ in lemma \ref{lem:s4t1}, then
\begin{align}
&\quad \biggr\|\sup_{0<t<2^{-ak/\min\{1,\gamma\}}}
  \frac{|P_{a,\gamma}^t (f_k)(\Gamma(x,t))-f_k(\Gamma(x,t))|}{t^{\delta\min\{1,\gamma\}/a}}\biggr\|_{L^q (B(x_0 ,R))} \notag\\
&\lesssim \sum_{j=1}^\infty \frac{2^{j/2}}{j!}2^{-akj+k\delta}\biggr\|\sup_{0<t<2^{-k/\beta}}
  \biggr|\int_{\mathbb R}e^{i\Gamma(x,t)\cdot \xi}|\xi|^{aj}\hat{f_k}(\xi)\rmd \xi\biggr|\biggr\|_{L^q (B(x_0 ,R))} \notag\\
&\lesssim \sum_{j=1}^\infty \frac{2^{j/2}}{j!}2^{-akj+k\delta}\sum_{l\in \mathbb Z}\frac{1}{(1+|l|)^2}
  \biggr\|\int_{\mathbb R}e^{i(x+l/2^{k})\cdot \xi}|\xi|^{aj}\hat{f_k}(\xi)\rmd \xi\biggr\|_{L^q (B(x_0 ,R))} \notag\\
&\lesssim \sum_{j=1}^\infty \frac{2^{j/2}}{j!}2^{k\delta}\|\hat{f_k}\|_{L^2 (\mathbb R)}
  \lesssim 2^{-s_1 k}\|f\|_{H^{s_1 +\delta}(\mathbb R)}.\label{eq:s4t14}
\end{align}
Inequalities (\ref{eq:s4t9})-(\ref{eq:s4t12}) and (\ref{eq:s4t14}) yield for $k\gg 1$,
\begin{equation}\label{eq:s4t15}
\biggr\|\sup_{0<t<1}\frac{|P_{a,\gamma}^t (f_k)(\Gamma(x,t))-f_k(x)|}{t^{\delta\min\{1,\gamma\}/a}}\biggr\|_{L^q (B(x_0 ,R))}
\lesssim 2^{-\varepsilon k/2}\|f\|_{H^{s_1+\delta}(\mathbb R)}.
\end{equation}
It is obvious that (\ref{eq:s4t1}) follows from (\ref{eq:s4t7}), (\ref{eq:s4t8}) and (\ref{eq:s4t15}).\\
\quad $(2)$ The proof is nearly identical to that of $(1)$. It suffices to show that for $q=\min\{p,2\}$
and any $\varepsilon >0$, $s_1=s_0+\varepsilon$, we have
\begin{equation}\label{eq:s4t16}
\biggr\|\sup_{0<t<1}\frac{|P_{a,\gamma}^t (f)(\Gamma(x,t))-f(x)|}{t^{\beta\delta}}\biggr\|_{L^q (B(x_0 ,R))}
\leq C_{\varepsilon}\|f\|_{H^{s_1+\delta}(\mathbb R)},
\end{equation}
where $0<\delta<\min\{\beta,\gamma\}/\beta$. We decompose $f$ as before and get
\begin{align}
&\quad \biggr\|\sup_{0<t<1}\frac{|P_{a,\gamma}^t (f)(\Gamma(x,t))-f(x)|}{t^{\beta\delta}}\biggr\|_{L^q (B(x_0 ,R))} \notag\\
&\leq \sum_{k=0}^\infty \biggr\|\sup_{0<t<1}\frac{|P_{a,\gamma}^t (f_k)(\Gamma(x,t))-f_k(x)|}{t^{\beta\delta}}\biggr\|_{L^q (B(x_0 ,R))}.\label{eq:s4t17}
\end{align}
For the case $k\lesssim 1$, by an argument similar to (\ref{eq:s4t8}), we obtain
\begin{equation}\label{eq:s4t18}
\biggr\|\sup_{0<t<1}\frac{|P_{a,\gamma}^t (f_k)(\Gamma(x,t))-f_k(x)|}{t^{\beta\delta}}\biggr\|_{L^q (B(x_0 ,R))}
 \lesssim\|f\|_{H^{s_1 +\delta}(\mathbb R)} .
\end{equation}
For the case $k\gg 1$, we have
\begin{align}
&\quad \biggr\|\sup_{0<t<1}\frac{|P_{a,\gamma}^t (f_k)(\Gamma(x,t))-f_k(x)|}{t^{\beta\delta}}\biggr\|_{L^q (B(x_0 ,R))} \notag\\
&\leq \biggr\|\sup_{2^{-k/\beta}\leq t<1}\frac{|P_{a,\gamma}^t (f_k)(\Gamma(x,t))-f_k(x)|}{t^{\beta\delta}}\biggr\|_{L^q (B(x_0 ,R))} \notag\\
&\qquad +\biggr\|\sup_{0<t<2^{-k/\beta}}\frac{|P_{a,\gamma}^t (f_k)(\Gamma(x,t))-f_k(x)|}{t^{\beta\delta}}\biggr\|_{L^q (B(x_0 ,R))} :=III+IV.\label{eq:s4t19}
\end{align}
We first estimate the part $III$. Note that $q=\min\{p,2\}$, we have
\begin{align}
III &\leq 2^{\delta k}\biggr\|\sup_{2^{-k/\beta}\leq t<1}|P_{a,\gamma}^t (f_k)(\Gamma(x,t))-f_k(x)|\biggr\|_{L^q (B(x_0 ,R))} \notag\\
&\lesssim 2^{\delta k} \biggr\{\biggr\| \sup_{2^{-k/\beta}\leq t<1}
  |P_{a,\gamma}^t (f_k)(\Gamma(x,t))| \biggr\|_{L^p (B(x_0 ,R))}+\|f_k\|_{L^2 (B(x_0,R))}\biggr\} \notag\\
&\lesssim 2^{(\delta +s_0+\varepsilon/2)k}\|f_k\|_{L^2 (\mathbb R)} \notag\\
&\lesssim 2^{-\varepsilon k/2}\|f\|_{H^{s_1+\delta}(\mathbb R)}.\label{eq:s4t20}
\end{align}
For the part $IV$, by Minkowski's inequality,
\begin{align}
IV &\leq \biggr\|\sup_{0<t<2^{-k/\beta}}\frac{|P_{a,\gamma}^t (f_k)(\Gamma(x,t))-f_k(\Gamma(x,t))|}{t^{\beta\delta}}\biggr\|_{L^q (B(x_0 ,R))} \notag\\
&\qquad +\biggr\|\sup_{0<t<2^{-k/\beta}}\frac{|f_k(\Gamma(x,t))-f_k(x)|}{t^{\beta\delta}}\biggr\|_{L^q (B(x_0 ,R))}.\label{eq:s4t21}
\end{align}
Arguing as in \cite[Theorem 1.2 (2)]{Li2021arxiv}, by Taylor's formula, we have
\begin{equation}\label{eq:s4t22}
\biggr\|\sup_{0<t<2^{-k/\beta}}\frac{|f_k(\Gamma(x,t))-f_k(x)|}{t^{\beta\delta}}\biggr\|_{L^q (B(x_0 ,R))}
 \lesssim 2^{-s_1 k}\|f\|_{H^{s_1+\delta}(\mathbb R)}.
\end{equation}
Whether $a\geq \min\{1,\gamma\}$ with $0<\beta<\min\{1,\gamma\}/a$,
or $a< \min\{1,\gamma\}$ with $0<\beta\leq 1$, by Lemma \ref{lem:s4t1} and (\ref{eq:s4t13}), we have
\begin{align}
&\quad \biggr\|\sup_{0<t<2^{-k/\beta}}
  \frac{|P_{a,\gamma}^t (f_k)(\Gamma(x,t))-f_k(\Gamma(x,t))|}{t^{\beta\delta}}\biggr\|_{L^q (B(x_0 ,R))} \notag\\
&\lesssim \sum_{j=1}^\infty \frac{2^{j/2}}{j!}2^{-\min\{1,\gamma\}kj/\beta}2^{k\delta}
  \biggr\|\sup_{0<t<2^{-k/\beta}}
   \biggr|\int_{\mathbb R}e^{i\Gamma(x,t)\cdot \xi}|\xi|^{aj}\hat{f_k}(\xi)\rmd \xi\biggr|\biggr\|_{L^q (B(x_0 ,R))} \notag\\
&\lesssim \sum_{j=1}^\infty \frac{2^{j/2}}{j!}2^{-\min\{1,\gamma\}kj/\beta}2^{k\delta}
  \sum_{l\in \mathbb Z}\frac{1}{(1+|l|)^2}
   \biggr\|\int_{\mathbb R}e^{i(x+l/2^{k})\cdot \xi}|\xi|^{aj}\hat{f_k}(\xi)\rmd \xi\biggr\|_{L^q (B(x_0 ,R))} \notag\\
&\lesssim \sum_{j=1}^\infty \frac{2^{j/2}}{j!}2^{k\delta}\|\hat{f_k}\|_{L^2 (\mathbb R)}
  \lesssim 2^{-s_1 k}\|f\|_{H^{s_1 +\delta}(\mathbb R)}.\label{eq:s4t23}
\end{align}
Inequalities (\ref{eq:s4t19})-(\ref{eq:s4t23}) yield for $k\gg 1$,
\begin{equation}\label{eq:s4t24}
\biggr\|\sup_{0<t<1}\frac{|P_{a,\gamma}^t (f_k)(\Gamma(x,t))-f_k(x)|}{t^{\beta\delta}}\biggr\|_{L^q (B(x_0 ,R))}
\lesssim 2^{-\varepsilon k/2}\|f\|_{H^{s_1+\delta}(\mathbb R)}.
\end{equation}
Then we arrive at (\ref{eq:s4t16}). \hfill $\Box$

\begin{Remark}
Let $\gamma >1$ and $a\geq 2$.
If the function $\Gamma$ satisfies (\ref{eq:s1t3}) for $1/a\leq \beta \leq 1$ and the bilipschitz condition, that is,
\[
  C_1 |x-y|\leq |\Gamma(x,t)-\Gamma(y,t)| \leq C_2 |x-y|
\]
uniformly for $x, y\in [-1,1]$ and $t\in [0,1]$,
then the regularity index $s_0$ in Theorem \ref{thm:s1t3} can be chosen as
$\min\{a(1-1/\gamma)/4,1/4\}$ according to {\cite[Theorem 1.3]{Niu2020curve}}.

We also note that Theorem \ref{thm:s1t3} still holds for higher dimensions $n\geq 2$, with almost the same proof.
\end{Remark}

\subsection{sharpness of Theorem \ref{thm:s1t4}}

Before we discuss the sharpness of Theorem \ref{thm:s1t4}, we first make a simple yet crucial observation.
From the proof of Theorem \ref{thm:s1t4}, it is evident that the following stronger version holds.

\begin{Theorem}\label{thm:s5t1}
Let $\gamma >0$, $a>0$ and $p\geq 1$. For any $\varepsilon >0$,
if a function $\hat{f}$ supported on $\{\xi:|\xi|\approx R\}$ for some $R \gg 1$ satisfies
\begin{equation}\label{eq:s5t1}
\biggr\|\sup_{0<t<1}|P_{a,\gamma}^t f|\biggr\|_{L^p (B(0,1))}
\leq C_{\varepsilon} R^{\varepsilon} \|f\|_{L^2 (\mathbb R)},
\end{equation}
where $C_{\varepsilon}>0$, then for any $0\leq \delta <a$, we have
\begin{equation}\label{eq:s5t2}
\biggr\|\sup_{0<t<1}\frac{|P_{a,\gamma}^t (f)(x)-f(x)|}{t^{\delta\min\{1,\gamma\}/a}}\biggr\|_{L^q (B(0,1))}
\leq C'_{\varepsilon} R^{\varepsilon+\delta} \|f\|_{L^2 (\mathbb R)}, \quad q=\min\{p,2\}.
\end{equation}
\end{Theorem}

The sharpness of Theorem \ref{thm:s1t4} is reflected in the following two aspects.
On the one hand, inspired by the work of Li and Wang \cite{Li2021arxiv},
we will give a counterexample to show that under the conditions of Theorem \ref{thm:s5t1},
$t^{\delta\min\{1,\gamma\}/a}$ in inequality (\ref{eq:s5t2}) can not be replaced by $t^{\delta w}$, where $w>\min\{1,\gamma\}/a$.
On the other hand, the following conclusion implies that for non-zero Schwartz functions, the convergence rate is no faster than
$t^{\min\{1,\gamma\}}$ as $t$ tends to zero along vertical lines.

\begin{Theorem}\label{thm:s5t2}
$(1)$ Let $a>0$, $\gamma>0$ and $0<\beta<1$. There exists
\[
  \Gamma(x,t)=x-t^{\beta}
\]
such that for each Schwartz function $f$, if
\begin{equation}\label{eq:s5t3}
\lim_{t\rightarrow 0^+}\frac{P_{a,\gamma}^t (f)(\Gamma(x,t))-f(x)}{t^{\min\{\beta,\gamma\}}}=0,\quad a.e.\quad x\in\mathbb R,
\end{equation}
then we have $f\equiv 0$.\\
\quad $(2)$ For each Schwartz function $f$, if
\begin{equation}\label{eq:s5t4}
\lim_{t\rightarrow 0^+}\frac{P_{a,\gamma}^t (f)(x)-f(x)}{t^{\min\{1,\gamma\}}}=0,\quad a.e.\quad x\in\mathbb R,
\end{equation}
then we have $f\equiv 0$.
\end{Theorem}

We are now in a position to prove the necessity of Theorem \ref{thm:s5t1}.

\proof[Necessity of Theorem \ref{thm:s5t1}.]
Let's consider the function $f_R$ defined as follows,
\[
  \widehat{f_R}(\xi)=\chi_{I_R}(\xi),
\]
where $I_R := [R,R+1]$. By H\"older's inequality, we have
\[
  \biggr\|\sup_{0<t<1}|P_{a,\gamma}^t f_R|\biggr\|_{L^p (B(0,1))}
   \lesssim_p \|f_R\|_{L^2 (\mathbb R)}, \quad \forall p\geq 1.
\]
Therefore, for any $\varepsilon >0$,
\[
  \biggr\|\sup_{0<t<1}|P_{a,\gamma}^t f_R|\biggr\|_{L^p (B(0,1))}
   \lesssim_p R^{\varepsilon} \|f_R\|_{L^2 (\mathbb R)}, \quad \forall p\geq 1.
\]
Fix $p\geq 1$ and suppose that there exist constants $\delta_1 \geq 0$ and $\delta_2 \geq 0$ satisfying
\begin{equation}\label{eq:s5t5}
\biggr\|\sup_{0<t<1}\frac{|P_{a,\gamma}^t (f_R)(x)-f_R (x)|}{t^{\delta_1}}\biggr\|_{L^q (B(0,1))}
\leq C'_{\varepsilon} R^{\varepsilon+\delta_2} \|f_R\|_{L^2 (\mathbb R)}, \quad q=\min\{p,2\}.
\end{equation}
Our goal is to show that $\delta_1 \leq \delta_2 \min\{1,\gamma\}/a$.
Consider the function $h(t)=t^2 +t^{2\gamma}$.
It is easy to see that $h$ is a strictly increasing function on $[0,\infty)$ which satisfies $h(0)=0$ and $h(1)=2$.
Therefore, there exists $t_0 \in (0,1)$ such that $h(t_0)=R^{-2a}/10^4$. Since
\[
  h\biggr(\frac{R^{-a/\min\{1,\gamma\}}}{100^{1/\min\{1,\gamma\}}}\biggr)
   =\frac{R^{-2a/\min\{1,\gamma\}}}{10^{4/\min\{1,\gamma\}}}
    +\frac{R^{-2a\gamma/\min\{1,\gamma\}}}{10^{4\gamma/\min\{1,\gamma\}}}>h(t_0),
\]
we obtain
\begin{equation}\label{eq:s5t6}
t_0 <\frac{R^{-a/\min\{1,\gamma\}}}{100^{1/\min\{1,\gamma\}}}.
\end{equation}
For almost every $x\in B(0,1/1000)$, by Taylor's formula,
\begin{align}
&\quad \frac{|P_{a,\gamma}^{t_0} (f_R)(x)-f_R (x)|}{{t_0}^{\delta_1}}
  \approx {t_0}^{-\delta_1}\biggr|\int_R^{R+1} (e^{it_0 |\xi|^a}e^{-t_0^{\gamma}|\xi|^a}-1)e^{ix\cdot \xi}\rmd \xi\biggr| \notag\\
&={t_0}^{-\delta_1}\biggr|\int_R^{R+1} e^{ix\cdot \xi}(it_0 |\xi|^a -t_0^{\gamma}|\xi|^a)\rmd \xi
  +\sum_{j\geq 2}\frac{1}{j!}\int_R^{R+1} e^{ix\cdot \xi}(it_0 |\xi|^a -t_0^{\gamma}|\xi|^a)^j\rmd \xi\biggr| \notag\\
&\geq {t_0}^{-\delta_1}\biggr|\biggr|\int_R^{R+1} e^{ix\cdot \xi}(it_0 |\xi|^a -t_0^{\gamma}|\xi|^a)\rmd \xi\biggr|
   -\biggr|\sum_{j\geq 2}\frac{1}{j!}\int_R^{R+1} e^{ix\cdot \xi}(it_0 |\xi|^a -t_0^{\gamma}|\xi|^a)^j\rmd \xi\biggr|\biggr|. \label{eq:s5t7}
\end{align}
Note that
\[
  \biggr|\int_R^{R+1} e^{ix\cdot \xi}(it_0 |\xi|^a -t_0^{\gamma}|\xi|^a)\rmd \xi\biggr|
   =h(t_0)^{1/2}\biggr|\int_0^1 e^{ix\cdot \eta} |\eta +R|^{a}\rmd \eta\biggr|.
\]
Thus, it follows from $|x \eta|<1/1000$ that for any $x\in B(0,1/1000)$,
\begin{align}
\biggr|\int_R^{R+1} e^{ix\cdot \xi}(it_0 |\xi|^a -t_0^{\gamma}|\xi|^a)\rmd \xi\biggr|
&\geq h(t_0)^{1/2} \int_0^1 (1-|x \eta|)|\eta +R|^{a}\rmd \eta \notag\\
&\geq \frac{h(t_0)^{1/2} R^a}{2}=\frac{1}{200}. \label{eq:s5t8}
\end{align}
Take $R$ sufficiently large that $(R+1)^a/R^a \leq 2$. Then
\begin{align}
\biggr|\sum_{j\geq 2}\frac{1}{j!}\int_R^{R+1} e^{ix\cdot \xi}(it_0 |\xi|^a -t_0^{\gamma}|\xi|^a)^j\rmd \xi\biggr|
&\leq \sum_{j\geq 2}\frac{1}{j!} h(t_0)^{j/2} \int_R^{R+1} |\xi|^{aj} \rmd \xi \notag\\
&\leq \sum_{j\geq 2}\frac{1}{100^j j!}\biggr( \frac{R+1}{R} \biggr)^{aj} \leq \frac{e-2}{2500}. \label{eq:s5t9}
\end{align}
Combining (\ref{eq:s5t6})-(\ref{eq:s5t9}), we get
\begin{align*}
&\quad \biggr\|\sup_{0<t<1}\frac{|P_{a,\gamma}^t (f_R)(x)-f_R (x)|}{t^{\delta_1}}\biggr\|_{L^q (B(0,1))} \notag\\
&\gtrsim \biggr\|\frac{|P_{a,\gamma}^{t_0} (f_R)(x)-f_R (x)|}{{t_0}^{\delta_1}}\biggr\|_{L^q (B(0,1/1000))}
  \gtrsim R^{a\delta_1 /\min\{1,\gamma\}}.
\end{align*}
Consequently, by inequality (\ref{eq:s5t5}),
\[
  R^{a\delta_1 /\min\{1,\gamma\}-\delta_2 -\varepsilon}
   \lesssim C'_{\varepsilon}\|f_R\|_{L^2 (\mathbb R)}\approx C'_{\varepsilon}.
\]
Fix $\varepsilon >0$. When $R$ tends to infinity, this inequality holds only if
\[
  a\delta_1 /\min\{1,\gamma\}-\delta_2 -\varepsilon \leq 0.
\]
By the arbitrariness of $\varepsilon$, we have $\delta_1 \leq \delta_2 \min\{1,\gamma\}/a$.
\endproof

\begin{Remark}
The proof of the scalar result is still applicable to higher dimensions $n\geq 2$
if we replace the interval $I_R$ by the cube $Q_R$ defined by
\[
  Q_R :=\{\xi:R\leq \xi_i \leq R+1,1\leq i\leq n\}.
\]
Indeed, we can easily find a positive number $t_0 < R^{-a/\min\{1,\gamma\}}/(100 n^{a/2})^{1/\min\{1,\gamma\}}$
satisfying $h(t_0)=R^{-2a}/(10^4 n^a)$, where the function $h$ is defined as above.
Then by a similar argument we see that the estimates (\ref{eq:s5t7})-(\ref{eq:s5t9}) still hold.
\end{Remark}

Finally, we give a proof of Theorem \ref{thm:s5t2}.
The original idea for the proof of Theorem \ref{thm:s5t2} originates from \cite[Theorem 5.2]{Li2021JFA}.\\
\textbf{Proof of Theorem \ref{thm:s5t2}.}
$(1)$ We may assume that $t\in (0,1)$. Then by Taylor's formula, we obtain
\begin{align}
&\quad |P_{a,\gamma}^t (f)(\Gamma(x,t))-f(x)| \notag\\
&\approx \biggr|\int_{\{\xi:t^\beta|\xi|+(t+t^{\gamma})|\xi|^a \leq \frac{1}{1000}\}}
  e^{ix\cdot \xi}(e^{-it^{\beta}\xi +it|\xi|^a}e^{-t^{\gamma}|\xi|^a}-1)\hat{f}(\xi)\rmd\xi \notag\\
&\quad + \int_{\{\xi:t^\beta|\xi|+(t+t^{\gamma})|\xi|^a > \frac{1}{1000}\}}
  e^{ix\cdot \xi}(e^{-it^{\beta}\xi +it|\xi|^a}e^{-t^{\gamma}|\xi|^a}-1)\hat{f}(\xi)\rmd\xi \biggr| \notag\\
&\geq \biggr|\biggr|\int_{\{\xi:t^\beta|\xi|+(t+t^{\gamma})|\xi|^a \leq \frac{1}{1000}\}}
  e^{ix\cdot \xi}(-it^{\beta}\xi +it|\xi|^a -t^{\gamma}|\xi|^a) \hat{f}(\xi)\rmd\xi \biggr| \notag\\
&\quad -\biggr|\sum_{j\geq 2}\frac{1}{j!}\int_{\{\xi:t^\beta|\xi|+(t+t^{\gamma})|\xi|^a \leq \frac{1}{1000}\}}
  e^{ix\cdot \xi}(-it^{\beta}\xi +it|\xi|^a -t^{\gamma}|\xi|^a)^j \hat{f}(\xi)\rmd\xi \notag\\
&\qquad + \int_{\{\xi:t^\beta|\xi|+(t+t^{\gamma})|\xi|^a > \frac{1}{1000}\}}
  e^{ix\cdot \xi}(e^{-it^{\beta}\xi +it|\xi|^a}e^{-t^{\gamma}|\xi|^a}-1)\hat{f}(\xi)\rmd\xi \biggr| \biggr|.\label{eq:s5t10}
\end{align}
It follows from $|e^{-it^{\beta}\xi +it|\xi|^a}e^{-t^{\gamma}|\xi|^a}-1|\lesssim 1$ that
\begin{align}
&\quad \biggr|\int_{\{\xi:t^\beta|\xi|+(t+t^{\gamma})|\xi|^a > \frac{1}{1000}\}}
  e^{ix\cdot \xi}(e^{-it^{\beta}\xi +it|\xi|^a}e^{-t^{\gamma}|\xi|^a}-1)\hat{f}(\xi)\rmd\xi \biggr| \notag\\
&\lesssim \int_{\mathbb R}(t^\beta|\xi|+t|\xi|^a +t^{\gamma}|\xi|^a)^2 |\hat{f}(\xi)|\rmd \xi. \label{eq:s5t11}
\end{align}
Furthermore, we have
\begin{align}
&\quad \biggr|\sum_{j\geq 2}\frac{1}{j!}\int_{\{\xi:t^\beta|\xi|+(t+t^{\gamma})|\xi|^a \leq \frac{1}{1000}\}}
  e^{ix\cdot \xi}(-it^{\beta}\xi +it|\xi|^a -t^{\gamma}|\xi|^a)^j \hat{f}(\xi)\rmd\xi \biggr| \notag\\
&\leq \sum_{j\geq 2}\frac{1}{j!}\int_{\{\xi:t^\beta|\xi|+(t+t^{\gamma})|\xi|^a \leq \frac{1}{1000}\}}
  (t^\beta|\xi|+t|\xi|^a +t^{\gamma}|\xi|^a)^2 |\hat{f}(\xi)|\rmd \xi \notag\\
&\lesssim \int_{\mathbb R}(t^\beta|\xi|+t|\xi|^a +t^{\gamma}|\xi|^a)^2 |\hat{f}(\xi)|\rmd \xi. \label{eq:s5t12}
\end{align}
Since $f$ is a Schwartz function, we have
\begin{align*}
&\quad \lim_{t\rightarrow 0^+}t^{-\min\{\beta,\gamma\}}\int_{\mathbb R}(t^\beta|\xi|+t|\xi|^a +t^{\gamma}|\xi|^a)^2 |\hat{f}(\xi)|\rmd \xi \notag\\
&\leq \lim_{t\rightarrow 0^+}t^{\min\{\beta,\gamma\}}\int_{\mathbb R}(|\xi|+2|\xi|^a)^2 |\hat{f}(\xi)|\rmd \xi =0.
\end{align*}
Based on (\ref{eq:s5t10})-(\ref{eq:s5t12}) and the Lebesgue dominated convergence theorem, we get
\begin{align*}
&\quad \lim_{t\rightarrow 0^+}\frac{|P_{a,\gamma}^t (f)(\Gamma(x,t))-f(x)|}{t^{\min\{\beta,\gamma\}}} \\
&\gtrsim \lim_{t\rightarrow 0^+}t^{-\min\{\beta,\gamma\}}\biggr|\int_{\{\xi:t^\beta|\xi|+(t+t^{\gamma})|\xi|^a \leq \frac{1}{1000}\}}
  e^{ix\cdot \xi}(-it^{\beta}\xi +it|\xi|^a -t^{\gamma}|\xi|^a) \hat{f}(\xi)\rmd\xi \biggr| \\
&=\biggr|\int_{\mathbb R}e^{ix\cdot \xi}Q(\xi) \hat{f}(\xi)\rmd\xi \biggr|,
\end{align*}
where
\begin{equation}
Q(\xi)=\begin{cases}
\xi,&\beta <\gamma; \\
|\xi|^a,&\beta >\gamma; \\
i\xi+|\xi|^a,&\beta =\gamma.
\end{cases}
\nonumber
\end{equation}
If $f$ is non-zero, then there exists a set of positive measure and a constant $C_1 >0$ such that for each $x$ in this set, it holds
\[
  \biggr|\int_{\mathbb R}e^{ix\cdot \xi}Q(\xi) \hat{f}(\xi)\rmd\xi \biggr|\geq C_1,
\]
which contradicts with (\ref{eq:s5t3}).

$(2)$ The proof is similar to $(1)$. By Taylor's formula, we get
\begin{align}
&\quad |P_{a,\gamma}^t (f)(x)-f(x)| \notag\\
&\approx \biggr|\int_{\mathbb R}e^{ix\cdot \xi}(e^{it|\xi|^a}e^{-t^{\gamma}|\xi|^a}-1)\hat{f}(\xi)\rmd \xi \biggr| \notag\\
&=\biggr|\int_{\{\xi:|it-t^{\gamma}|\cdot |\xi|^a \leq \frac{1}{1000}\}}
  e^{ix\cdot \xi}(e^{it|\xi|^a}e^{-t^{\gamma}|\xi|^a}-1)\hat{f}(\xi)\rmd\xi \notag\\
&\quad +\int_{\{\xi:|it-t^{\gamma}|\cdot|\xi|^a >\frac{1}{1000}\}}
  e^{ix\cdot \xi}(e^{it|\xi|^a}e^{-t^{\gamma}|\xi|^a}-1)\hat{f}(\xi)\rmd\xi \biggr| \notag\\
&\geq \biggr|\biggr|\int_{\{\xi:|it-t^{\gamma}|\cdot |\xi|^a \leq \frac{1}{1000}\}}
  e^{ix\cdot \xi}(it|\xi|^a -t^{\gamma}|\xi|^a)\hat{f}(\xi)\rmd\xi\biggr| \notag\\
&\quad -\biggr|\sum_{j\geq 2}\frac{1}{j!}\int_{\{\xi:|it-t^{\gamma}|\cdot |\xi|^a \leq \frac{1}{1000}\}}
  e^{ix\cdot \xi}(it|\xi|^a -t^{\gamma}|\xi|^a)^j \hat{f}(\xi)\rmd\xi \notag\\
&\qquad + \int_{\{\xi:|it-t^{\gamma}|\cdot |\xi|^a > \frac{1}{1000}\}}
  e^{ix\cdot \xi}(e^{it|\xi|^a}e^{-t^{\gamma}|\xi|^a}-1)\hat{f}(\xi)\rmd\xi \biggr| \biggr|.\label{eq:s5t13}
\end{align}
It is easy to see that
\begin{align}
&\quad \biggr|\int_{\{\xi:|it-t^{\gamma}|\cdot |\xi|^a > \frac{1}{1000}\}}
  e^{ix\cdot \xi}(e^{it|\xi|^a}e^{-t^{\gamma}|\xi|^a}-1)\hat{f}(\xi)\rmd\xi \biggr| \notag\\
&\lesssim (t^2 +t^{2\gamma})\int_{\mathbb R}|\xi|^{2a}|\hat{f}(\xi)|\rmd \xi.\label{eq:s5t14}
\end{align}
Furthermore, we also have
\begin{align}
&\quad \biggr|\sum_{j\geq 2}\frac{1}{j!}\int_{\{\xi:|it-t^{\gamma}|\cdot |\xi|^a \leq \frac{1}{1000}\}}
  e^{ix\cdot \xi}(it|\xi|^a -t^{\gamma}|\xi|^a)^j \hat{f}(\xi)\rmd\xi \biggr| \notag\\
&\leq \sum_{j\geq 2}\frac{1}{j!}\int_{\{\xi:|it-t^{\gamma}|\cdot |\xi|^a \leq \frac{1}{1000}\}}(t^2 +t^{2\gamma})|\xi|^{2a} |\hat{f}(\xi)|\rmd\xi\notag\\
&\lesssim (t^2 +t^{2\gamma})\int_{\mathbb R}|\xi|^{2a} |\hat{f}(\xi)|\rmd\xi .\label{eq:s5t15}
\end{align}
Combining (\ref{eq:s5t13})-(\ref{eq:s5t15}), we get
\begin{align*}
&\quad \lim_{t\rightarrow 0^+}\frac{P_{a,\gamma}^t (f)(x)-f(x)}{t^{\min\{1,\gamma\}}} \notag\\
&\gtrsim \lim_{t\rightarrow 0^+}\biggr|\int_{\{\xi:|it-t^{\gamma}|\cdot |\xi|^a \leq \frac{1}{1000}\}}
  e^{ix\cdot \xi}(it^{1-\min\{1,\gamma\}}|\xi|^a -t^{\gamma-\min\{1,\gamma\}}|\xi|^a)\hat{f}(\xi)\rmd\xi\biggr| \notag\\
&\approx \biggr| \int_{\mathbb R}e^{ix\cdot \xi}|\xi|^a \hat{f}(\xi) \biggr|.
\end{align*}
If $f$ is non-zero, then there exists a set of positive measure and a constant $C_2 >0$ such that for each $x$ in this set, it holds
\[
  \biggr|\int_{\mathbb R}e^{ix\cdot \xi} |\xi|^a \hat{f}(\xi)\rmd\xi \biggr|\geq C_2,
\]
which contradicts with (\ref{eq:s5t4}). \hfill $\Box$

\begin{Remark}
For higher dimensions, if we redefine the curve $\Gamma (x,t)$ as
\[
  \Gamma(x,t)=x-e_1 \cdot t^{\beta}
\]
where $e_1=(1,0,\cdots,0)$, then we get the same upper bound for the convergence rate as in Theorem \ref{thm:s5t2}.
Unlike the scalar model, we need to decompose the integral region in (\ref{eq:s5t10}) into the union of
$\Omega_1 :=\{\xi\in \mathbb R^n :t^\beta|\xi_1|+(t+t^{\gamma})|\xi|^a \leq 1/1000\}$ and
$\Omega_2 :=\{\xi\in \mathbb R^n :t^\beta|\xi_1|+(t+t^{\gamma})|\xi|^a > 1/1000\}$.
\end{Remark}

\appendix
\section{proof of Theorem \ref{thm:s1t4}.}
In this section, we give a proof of Theorem \ref{thm:s1t4}.
\proof By a standard argument, it suffices to show that for $q=\min\{p,2\}$
and any $x_0 \in \mathbb R$, $\varepsilon >0$, $s_1=s_0+\varepsilon$, we have
\begin{equation}\label{eq:s6t1}
\biggr\|\sup_{0<t<1}\frac{|P_{a,\gamma}^t (f)(x)-f(x)|}{t^{\delta\min\{1,\gamma\}/a}}\biggr\|_{L^q (B(x_0 ,1))}
\leq C_{\varepsilon}\|f\|_{H^{s_1+\delta}(\mathbb R)}.
\end{equation}
By translation, the problem is reduced to showing that
\begin{equation}\label{eq:s6t2}
\biggr\|\sup_{0<t<1}\frac{|P_{a,\gamma}^t (f)(x)-f(x)|}{t^{\delta\min\{1,\gamma\}/a}}\biggr\|_{L^q (B(0,1))}
\leq C_{\varepsilon}\|f\|_{H^{s_1+\delta}(\mathbb R)}.
\end{equation}
Indeed, if (\ref{eq:s6t2}) holds for all $f\in H^{s_1+\delta}(\mathbb R)$, then take $f_0$ defined by
\[
  \hat{f_0}(\xi)=e^{ix_0 \cdot \xi}\hat{f}(\xi)
\]
and replace $f$ by $f_0$ in (\ref{eq:s6t2}). It is obvious that (\ref{eq:s6t1}) follows from changing of variables.
To prove (\ref{eq:s6t2}), we decompose $f$ as follows:
\[
  f=\sum_{k=0}^\infty f_k,
\]
where $supp\hat{f_0}\subset B(0,1)$, $supp \hat{f_k}\subset \{\xi:|\xi|\approx 2^k\}$, $k\geq 1$. Therefore,
\begin{equation}\label{eq:s6t3}
\biggr\|\sup_{0<t<1}\frac{|P_{a,\gamma}^t (f)(x)-f(x)|}{t^{\delta\min\{1,\gamma\}/a}}\biggr\|_{L^q (B(0,1))}
\leq \sum_{k=0}^\infty \biggr\|\sup_{0<t<1}\frac{|P_{a,\gamma}^t (f_k)(x)-f_k (x)|}{t^{\delta\min\{1,\gamma\}/a}}\biggr\|_{L^q (B(0,1))}.
\end{equation}
For the case $k\lesssim 1$, similar to (\ref{eq:s4t8}), we get
\begin{equation}\label{eq:s6t4}
\biggr\|\sup_{0<t<1}\frac{|P_{a,\gamma}^t(f_k)(x)-f_k (x)|}{t^{\delta\min\{1,\gamma\}/a}}\biggr\|_{L^q (B(0,1))}
\lesssim\|f\|_{H^{s_1+\delta}(\mathbb R)}.
\end{equation}
For the case $k\gg 1$, we have
\begin{align}
&\quad \biggr\|\sup_{0<t<1}\frac{|P_{a,\gamma}^t(f_k)(x)-f_k (x)|}{t^{\delta\min\{1,\gamma\}/a}}\biggr\|_{L^q (B(0,1))} \notag\\
&\leq \biggr\|\sup_{2^{-ak/\min\{1,\gamma\}}\leq t<1}\frac{|P_{a,\gamma}^t(f_k)(x)-f_k (x)|}{t^{\delta\min\{1,\gamma\}/a}}\biggr\|_{L^q (B(0,1))} \notag\\
&\qquad +\biggr\|\sup_{0<t<2^{-ak/\min\{1,\gamma\}}}\frac{|P_{a,\gamma}^t(f_k)(x)-f_k (x)|}{t^{\delta\min\{1,\gamma\}/a}}\biggr\|_{L^q (B(0,1))}:= V+VI. \label{eq:s6t5}
\end{align}
For the part $V$, similar to (\ref{eq:s4t10}), we can easily get
\begin{equation}\label{eq:s6t6}
V \lesssim 2^{-\varepsilon k/2} \|f\|_{H^{s_1 +\delta}(\mathbb R)}.
\end{equation}
Finally, we estimate the part $VI$. By (\ref{eq:s4t13}), we obtain
\[
  |P_{a,\gamma}^t (f_k)(x)-f_k(x)|\lesssim \sum_{j=1}^\infty \frac{2^{j/2}}{j!}t^{\min\{1,\gamma\}j}
  \biggr|\int_{\mathbb R}e^{ix\cdot \xi}|\xi|^{aj}\hat{f_k}(\xi)\rmd \xi\biggr|,
\]
which implies that
\begin{align}
VI &=\biggr\|\sup_{0<t<2^{-ak/\min\{1,\gamma\}}}\frac{|P_{a,\gamma}^t (f_k)(x)-f_k(x)|}{t^{\delta\min\{1,\gamma\}/a}}\biggr\|_{L^q (B(0,1))} \notag\\
&\lesssim \sum_{j=1}^\infty \frac{2^{j/2}}{j!} 2^{\delta k}2^{-akj}
  \biggr\|\int_{\mathbb R}e^{ix\cdot \xi}|\xi|^{aj}\hat{f_k}(\xi)\rmd \xi\biggr\|_{L^2 (B(0,1))} \notag\\
&\lesssim \sum_{j=1}^\infty \frac{2^{j/2}}{j!} 2^{\delta k} \|\hat{f_k}\|_{L^2 (\mathbb R)}
  \lesssim 2^{-s_1 k} \|f\|_{H^{s_1 +\delta}(\mathbb R)}.\label{eq:s6t7}
\end{align}
It follows from (\ref{eq:s6t5})-(\ref{eq:s6t7}) that for the case $k\gg 1$, we have
\[
  \biggr\|\sup_{0<t<1}\frac{|P_{a,\gamma}^t(f_k)(x)-f_k (x)|}{t^{\delta\min\{1,\gamma\}/a}}\biggr\|_{L^q (B(0,1))}
  \lesssim 2^{-\varepsilon k/2} \|f\|_{H^{s_1 +\delta}(\mathbb R)}.
\]
Then we arrive at (\ref{eq:s6t2}).
\endproof


\end{document}